\documentclass[11pt,reqno]{amsart}
\usepackage[top=26mm, bottom=26mm, left=23mm, right=23mm]{geometry}
\usepackage{graphicx}
\usepackage{color}
\usepackage{multirow}
\newcommand{\lsp}{\vspace{3mm}}

\usepackage{times}

\newcommand{\vtwo}[2]{\left[\begin{array}{c} #1 \\ #2 \end{array}\right]}
\newcommand{\vfour}[4]{\left[\begin{array}{c} #1 \\ #2 \\ #3 \\ #4 \end{array}\right]}

\usepackage{bm}

\newcommand{\vct}[1]{\bm{\mathsf{#1}}}
\newcommand{\pvct}[1]{\bm{#1}}
\newcommand{\mtx}[1]{\bm{\mathsf{#1}}}
\newcommand{\mtb}[1]{\bm{\mathcal{#1}}}

\newtheorem{remark}{Remark}
\newtheorem{definition}{Definition}

\begin{document}

\begin{center}
\textbf{\large{Compressing rank-structured matrices via randomized sampling}}

\vspace{2mm}

\textit{Per-Gunnar Martinsson, Dept.~of Applied Math., Univ.~of Colorado Boulder}

\vspace{2mm}

March 22, 2015

\vspace{2mm}

\begin{minipage}{127mm}
\textbf{Abstract:} Randomized sampling has recently been
proven a highly efficient technique for computing
approximate factorizations of matrices that have low numerical
rank. This paper describes an extension of such techniques to a
wider class of matrices that are not themselves rank-deficient,
but have off-diagonal blocks that are; specifically, the classes
of so called \textit{Hierarchically Off-Diagonal Low Rank (HODLR)}
matrices and \textit{Hierarchically Block Separable (HBS)} matrices
(a.k.a.~``Hierarchically Semi-Separable (HSS)'' matrices).
Such matrices arise
frequently in numerical analysis and signal processing, in particular
in the construction of fast methods for solving differential and
integral equations numerically. These structures admit algebraic operations
(matrix-vector multiplications, matrix factorizations, matrix inversion,
etc.) to be performed very rapidly; but only once a data-sparse
representation of the matrix has been constructed. How to rapidly
compute this representation in the first place is much less well
understood. The present paper demonstrates that if an $N\times N$
matrix can be applied to a vector in $O(N)$ time, and if the ranks
of the off-diagonal blocks are bounded by an integer $k$, then the
cost for constructing a HODLR representation is $O(k^{2}\,N\,(\log N)^{2})$,
and the cost for constructing an HBS representation is $O(k^{2}\,N\,\log N)$
(assuming of course, that the matrix is compressible in the respective format).
The point is that when legacy codes (based on, e.g., the Fast Multipole Method)
can be used for the fast matrix-vector multiply, the proposed algorithm can
be used to obtain the data-sparse representation of the matrix, and then well-established
techniques for HODLR/HBS matrices can be used to invert or factor the matrix.
The proposed scheme is also useful in simplifying the implementation of certain
operations on rank-structured matrices such as the matrix-matrix multiplication,
low-rank update, addition, etc.
%
\end{minipage}

\end{center}

\section{Introduction}

A ubiquitous task in computational science is to rapidly perform
linear algebraic operations involving very large matrices. Such
operations typically exploit special ``structure'' in the matrix
since the costs of standard techniques tend to scale prohibitively
fast with matrix size; for a general $N\times N$ matrix, it costs
$O(N^{2})$ operations to perform a matrix-vector multiplication,
$O(N^{3})$ operations to perform Gaussian elimination or to invert
the matrix, etc. A well-known form of ``structure'' in a
matrix is sparsity. When at most a few entries in each row of the
matrix are non-zero (as is the case, e.g., for matrices
arising upon the discretization of differential equations, or
representing the link structure of the World Wide Web) matrix-vector
multiplications can be performed in $O(N)$ operations instead of
$O(N^{2})$. The description ``data-sparse'' applies to a matrix that
may be dense, but that shares the key characteristic of a sparse
matrix that some linear algebraic operations, typically the
matrix-vector multiplication, can to high precision be executed in
fewer than $O(N^{2})$ operations (often in close to linear time).

Several different formats for rank-structured matrices have been
proposed in the literature. In this manuscript, we rely on the
so called \textit{Hierarchically Off-Diagonal Low Rank (HODLR)}
format. This name was minted in \cite{2013_darve_FDS}, but this
class of matrices has a long history. It is a special case of
the $\mathcal{H}$-matrix format introduced by Hackbusch and
co-workers \cite{hackbusch,2008_bebendorf_book}, and was used
explicitly in \cite[Sec.~4]{2009_martinsson_FEM}. The HODLR format
is very easy to describe and easy to use, but can lead to less than
optimal performance due to the fact that the basis matrices used
to represent large blocks are stored explicitly, leading to a
$O(k\,N\,\log(N))$ storage requirement for a HODLR matrix whose
off-diagonal blocks have rank at most $k$. To attain \textit{linear}
storage requirements and arithmetic operations, one can switch to
a format that expresses all basis matrices \textit{hierarchically};
in other words, the basis matrices used on one level are expressed
implicitly in terms of the basis matrices on the next finer level.
We sometimes say that we use \textit{nested} basis matrices. To
be precise, we use the \textit{Hierarchically Block Separable (HBS)}
format that was described in \cite{2005_martinsson_fastdirect,2012_martinsson_FDS_survey}.
This format is closely related to the \textit{Hierarchically Semi-Separable (HSS)}
\cite{2005_gu_HSS,2009_xia_superfast} format, and is also related to
the $\mathcal{H}^{2}$-matrix format \cite{2002_hackbusch_H2,2010_borm_book}.

The most straight-forward technique for computing a data-sparse representation
of a rank-structured $N\times N$ matrix $\mtx{A}$ is to explicitly form all matrix
elements, and then to compress the off-diagonal blocks using, e.g., the SVD.
This approach can be executed stably \cite{2010_gu_xia_HSS,2002_hackbusch_borm_datasparse},
but it is often prohibitively expensive, with an $O(k\,N^{2})$ asymptotic cost, where $k$ is the
rank of the off-diagonal blocks (in the HSS-sense).
Fortunately, there exist for specific applications much faster methods for constructing
HSS representations.
When the matrix $\mtx{A}$ approximates a boundary integral operator in the plane,
the technique of \cite{2005_martinsson_fastdirect} computes a representation in $O(k^{2}\,N)$
time by exploiting representation results from potential theory.
In other environments, it is possible to use known regularity properties of
the off-diagonal blocks in conjunction with interpolation techniques
to obtain rough initial factorizations, and then recompress these to
obtain factorizations with close to optimal ranks \cite{2010_borm_book,2007_martinsson_hudson}.
A particularly popular version of the ``regularity + recompression'' method
is the so called \textit{Adaptive Cross Approximation} technique which was
initially proposed for $\mathcal{H}$-matrices
\cite{2000_bebendorf_ACA,2005_borm_grasedyck_ACA,huygens} but has recently been
modified to obtain a representation of a matrix in a format similar to the HSS
\cite{2010_frederix}.


The purpose of the present paper is to describe a fast and simple
randomized technique for computing a data sparse representation of
a rank-structured matrix which can rapidly be applied to a vector.
The existence of such a technique means that the advantages of the
HODLR and HBS formats --- fast inversion
and factorization algorithms
in particular --- become available for any matrix that can currently
be applied via the FMM, via an $\mathcal{H}$-matrix calculation, or
by any other existing data-sparse format (provided of course that
the matrix is in principle rank-structured). In order
to describe the cost of the algorithm precisely, we introduce some
notation: We let $\mtx{A}$ be an $N\times N$ matrix whose off-diagonal
blocks have maximal rank $k$, we let $T_{\rm mult}$ denote the time required to
perform a matrix-vector multiplication $\vct{x} \mapsto \mtx{A}\,\vct{x}$ or
$\mtx{x}\mapsto \mtx{A}^{*}\,\vct{x}$, we let $T_{\rm rand}$ denote the cost of
constructing a pseudo random number from a normalized Gaussian
distribution, and $T_{\rm flop}$ denote
the cost of a floating point operation. The computational cost
$T_{\rm total}$ of the algorithm for the HBS format then satisfies
\begin{equation}
\label{eq:Ttotal}
T_{\rm total} \sim T_{\rm mult}\times k\,\log(N) +
T_{\rm rand} \times (k+p)\,N\,\log(N) +
T_{\rm flop}\times k^{2}\,N\,\log(N),
\end{equation}
where $p$ is a tuning parameter
that balances computational cost against the probability of not
meeting the requested accuracy. Setting $p=10$ is often a good
choice which leads to a ``failure probability'' of less that $10^{-9}$,
see Remark \ref{remark:tuning}. In particular, if $T_{\rm mult}$ is
$O(N)$, then the method presented here has overall complexity $O(k^{2}\,N\,\log(N))$.
For the HODLR format, an additional factor of $\log N$ arises, cf.~(\ref{eq:T_HODLR}).

The work presented is directly inspired by \cite{2011_martinsson_randomhudson}
(which is based on a 2008 preprint \cite{2008_martinsson_randomhudson}), which
presented a similar algorithm with $O(k^{2}\,N)$ complexity for the compression
of an HBS matrix. This is better by a factor of $\log(N)$ compared to the present
work, but the algorithm of \cite{2011_martinsson_randomhudson} has a serious
limitation in that it requires the ability to evaluate $O(k\,N)$ entries of the
matrix to be compressed. This algorithm was later refined by Xia \cite{2013_xia_randomized}
and applied to the task of accelerated the ``nested dissection'' direct solver
\cite{george_1973,1989_directbook_duff} for elliptic PDEs to $O(N)$ complexity.
In 2011, by L.~Lin, J.~Lu, and L.~Ying presented an alternative algorithm
\cite{2011_lin_lu_ying} that
interacts with the matrix \textit{only} via the matrix-vector multiplication,
which makes the randomized compression idea much more broadly applicable than
the algorithm in \cite{2011_martinsson_randomhudson}. However, this came at the
cost of requiring $O(k\,\log(N))$ matrix-vector multiplies (just like the present
work). The algorithm proposed here is an evolution of \cite{2011_lin_lu_ying}
that does away with a step of least-square fitting that led to a magnification
of the sampling error resulting from the randomized approximation. Moreover,
we here present a new strategy for enforcing the ``nestedness'' of the basis
matrices that is required in the HBS format.

\begin{remark}
\label{remark:tuning}
The technique described in this paper utilizes a method for computing approximate
low-rank factorizations of matrices that is based on randomized sampling
\cite{2006_martinsson_random1_orig,2011_martinsson_random1,2011_martinsson_randomsurvey}. As a consequence, there is in
principle a non-zero risk that the method may fail to produce full accuracy in
any given realization of the algorithm. This risk can be controlled by
the user via the choice of the tuning parameter $p$ in (\ref{eq:Ttotal}),
for details see Section \ref{sec:randcomp}. Moreover, unlike
some better known randomized algorithms such as Monte Carlo, the accuracy
of the output of the algorithms under discussion here is typically very high;
in the environment described in the present paper, approximation errors of less
than $10^{-10}$ are entirely typical.
\end{remark}

%
%
%

\section{Preliminaries}

\subsection{Notation}
Throughout the paper, we measure vectors in $\mathbb{R}^{n}$ using their
Euclidean norm. The default norm for matrices will be the corresponding
operator norm $\|\mtx{A}\| = \sup_{\|\vct{x}\|=1}\|\mtx{A}\vct{x}\|$,
although we will sometimes also use the Frobenius norm
$\|\mtx{A}\|_{\rm Fro} = \left(\sum_{i,j}|\mtx{A}(i,j)|^{2}\right)^{1/2}$.

We use the notation of Golub and Van Loan \cite{golub} to specify
submatrices. In other words, if $\mtx{B}$ is an $m\times n$ matrix with
entries $b_{ij}$, and $I = [i_{1},\,i_{2},\,\dots,\,i_{k}]$ and $J =
[j_{1},\,j_{2},\,\dots,\,j_{\ell}]$ are two index vectors, then we
let $\mtx{B}(I,J)$ denote the $k\times \ell$ matrix
$$
\mtx{B}(I,J) = \left[\begin{array}{cccc}
b_{i_{1}j_{1}} & b_{i_{1}j_{2}} & \cdots & b_{i_{1}j_{\ell}} \\
b_{i_{2}j_{1}} & b_{i_{2}j_{2}} & \cdots & b_{i_{2}j_{\ell}} \\
\vdots         & \vdots         &        & \vdots         \\
b_{i_{k}j_{1}} & b_{i_{k}j_{2}} & \cdots & b_{i_{k}j_{\ell}}
\end{array}\right].
$$
We let $\mtx{B}(I,:)$ denote the matrix $\mtx{B}(I,[1,\,2,\,\dots,\,n])$,
and define $\mtx{B}(:,J)$ analogously.

The transpose of $\mtx{B}$ is denoted $\mtx{B}^{*}$, and we say that a matrix $\mtx{U}$
is \textit{orthonormal} if its  columns form an orthonormal set, so that $\mtx{U}^{*}\mtx{U} = I$.

\subsection{The QR factorization}
\label{sec:QR}
Any $m\times n$ matrix $\mtx{A}$ admits a \textit{QR factorization} of the form
\begin{equation}
\label{eq:QR}
\begin{array}{ccccccccccccccccccccc}
\mtx{A} & \mtx{P} &=& \mtx{Q} & \mtx{R},\\
m\times n & n \times n && m\times r & r \times n
\end{array}
\end{equation}
where $r = \min(m,n)$, $\mtx{Q}$ is orthonormal, $\mtx{R}$ is upper triangular,
and $P$ is a permutation matrix. The permutation matrix $\mtx{P}$ can more efficiently
be represented via a vector $J_c \in \mathbb{Z}_{+}^{n}$ of column indices
such that $\mtx{P} = \mtx{I}(:,J_c)$ where $\mtx{I}$ is the $n\times n$ identity matrix. As a
result, the factorization can be written as:
\begin{equation*}
\begin{array}{ccccccccccccccccccccc}
\mtx{A}(\colon,J_c) &=& \mtx{Q} & \mtx{R},\\
m\times n && m\times r & r \times n
\end{array}
\end{equation*}
The QR-factorization
is often built incrementally via a greedy algorithm such as column pivoted Gram-Schmidt.
This allows one to stop after the first $k$ terms have been computed to obtain a
``partial QR-factorization of $\mtx{A}$'':
\begin{equation*}
\begin{array}{ccccccccccccccccccccc}
\mtx{A}(\colon,J_c) &=& \left[\begin{matrix} \mtx{Q}^{(1)} & \mtx{Q}^{(2)} \end{matrix}\right] & \left[\begin{matrix} \mtx{R}^{(1)} \\ \mtx{R}^{(2)} \end{matrix}\right],\\
m\times n && m\times r & r \times n
\end{array}
\end{equation*}
That is, taking the first $k$ columns of $\mtx{Q}$ and the first $k$ rows of
$\mtx{R}$, we can obtain the approximation:
\begin{equation}
\label{eq:QR3}
\mtx{A}(\colon,J_c) \approx \mtx{Q}_k \mtx{R}_k
\end{equation}
We note that the partial factors $\mtx{Q}_k$ and $\mtx{R}_k$  can be obtained after
$k$ steps of the pivoted QR algorithm, without having to compute the
full matrices $\mtx{Q}$ and $\mtx{R}$.

\subsection{The singular value decomposition (SVD)}
\label{sec:SVD}
Let $\mtx{A}$ denote an $m\times n$ matrix, and set $r = \min(m,n)$.
Then $\mtx{A}$ admits a factorization
\begin{equation}
\label{eq:svd}
\begin{array}{ccccccccccccccccccccc}
\mtx{A} &=& \mtx{U} & \mtx{\Sigma} & \mtx{V}^{*},\\
m\times n && m\times r & r \times r & r\times n
\end{array}
\end{equation}
where the matrices $\mtx{U}$ and $\mtx{V}$ are orthonormal, and $\mtx{\Sigma}$ is diagonal.
We let $\{\vct{u}_{i}\}_{i=1}^{r}$ and $\{\vct{v}_{i}\}_{i=1}^{r}$ denote the columns of
$\mtx{U}$ and $\mtx{V}$, respectively. These vectors are the left and right singular vectors
of $\mtx{A}$. The diagonal elements $\{\sigma_{j}\}_{j=1}^{r}$ of
$\mtx{\Sigma}$ are the singular values of $\mtx{A}$. We order these so that
$\sigma_{1}  \geq \sigma_{2} \geq \cdots \geq \sigma_{r} \geq 0$.
We let $\mtx{A}_{k}$ denote the truncation of the SVD to its first $k$ terms, so that
$\mtx{A}_{k} = \sum_{j=1}^{k}\sigma_{j}\,\vct{u}_{j}\,\vct{v}_{j}^{*}$.
It is easily verified that
$$
\|\mtx{A} - \mtx{A}_{k}\|_{\rm spectral} = \sigma_{k+1},
\qquad\mbox{and that}\qquad
\|\mtx{A} - \mtx{A}_{k}\|_{\rm Fro} = \left(\sum_{j=k+1}^{\min(m,n)} \sigma_{j}^{2}\right)^{1/2},
$$
where $\|\cdot\|_{\rm spectral}$ denotes the operator norm and
$\|\cdot\|_{\rm Fro}$ denotes the Frobenius norm. Moreover,
the Eckart-Young theorem states that these errors are the smallest possible
errors that can be incurred when approximating $\mtx{A}$ by a matrix of rank $k$.


\subsection{The interpolatory decomposition (ID)}
\label{sec:ID}
Let $\mtx{A}$ denote an $m\times n$ matrix of rank $k$.
Then $\mtx{A}$ admits the factorization
$$
\begin{array}{ccccccccccccc}
\mtx{A} &=& \mtx{A}(:,J) & \mtx{X},\\
m\times n && m\times k & k\times n
\end{array}
$$
where $J$ is a vector of indices marking $k$ of the columns of $\mtx{A}$,
and the $k\times n$ matrix $\mtx{X}$ has the $k\times k$ identity matrix
as a submatrix and has the property that all its entries are
bounded by $1$ in magnitude. In other words, the interpolative
decomposition picks $k$ columns of $\mtx{A}$ as a basis for the column
space of $\mtx{A}$ and expresses the remaining columns in terms of the
chosen ones. The ID can be viewed as a modification to so the called
\textit{Rank-Revealing QR factorization} \cite{1987_chan_rrqr}. It can be computed in a stable
and accurate manner using the techniques of \cite{gu1996}, as described
in \cite{2005_martinsson_skel}. (Practical algorithms for
computing the interpolative decomposition produce a matrix $\mtx{X}$
whose elements slightly exceed $1$ in magnitude.)

\subsection{Randomized compression}
\label{sec:randcomp}
Let $\mtx{A}$ be a given $m\times n$ matrix that
can accurately be approximated by a matrix of rank $k$,
and suppose that we seek to determine a
matrix $\mtx{Q}$ with orthonormal columns (as few as possible) such that
$$
||\mtx{A} - \mtx{Q}\,\mtx{Q}^{*}\,\mtx{A}||
$$
is small. In other words, we seek a matrix $\mtx{Q}$ whose columns form an approximate orthornomal
basis (ON-basis) for the column space of $\mtx{A}$.
This task can efficiently be solved via the following randomized procedure:
\begin{enumerate}
\item Pick a small integer $p$ representing how much ``over-sampling'' we do.
      ($p=10$ is often good.)
\item Form an $n\times (k+p)$ matrix $\mtx{\Omega}$ whose entries are iid normalized Gaussian random numbers.
\item Form the product $\mtx{Y} = \mtx{A}\,\mtx{\Omega}$.
\item Construct a matrix $\mtx{Q}$ whose columns form an ON-basis for the columns of $\mtx{Y}$.
\end{enumerate}
Note that each column of the ``sample'' matrix $\mtx{Y}$ is a random
linear combination of the columns of $\mtx{A}$. We would therefore expect
the algorithm described to have a high probability of producing an
accurate result when $p$ is a large number. It is perhaps less
obvious that this probability depends only on $p$ (not on $m$ or $n$, or
any other properties of $\mtx{A}$), and that it approaches $1$ extremely
rapidly as $p$ increases. In fact, one can show that the basis $\mtx{Q}$
determined by the scheme above satisfies
\begin{equation}
\label{eq:err_prob}
\| \mtx{A} - \mtx{Q}\,\mtx{Q}^{*}\,\mtx{A} \|
    \leq \left[ 1 + 11 \sqrt{k+p} \cdot\sqrt{\min\{m,n\}} \right] \sigma_{k+1},
\end{equation}
with probability at least $1 - 6\cdot p^{-p}$, see \cite[Sec.~1.5]{2011_martinsson_randomsurvey}.
The error bound (\ref{eq:err_prob}) indicates that the error
produced by the randomized sampling procedure can be larger
than the theoretically minimal error $\sigma_{k+1}$ by
a factor of $1 + 11 \sqrt{k+p} \cdot\sqrt{\min\{m,n\}}$.
This crude bound is typically very pessimistic; for
specific situations sharper bounds have been proved, see \cite{2011_martinsson_randomsurvey}.

\begin{definition}
Let $\mtx{A}$ be an $m\times n$ matrix, and let $\varepsilon$ be a positive
number. We say that an $m\times \ell$ matrix $\mtx{Y}$
an $\varepsilon$-\textit{spanning matrix} for $\mtx{A}$, if
$\|\mtx{A} - \mtx{Y}\mtx{Y}^{\dagger}\mtx{A}\| \leq \varepsilon$.
Informally, this means that the columns of $\mtx{Y}$ span the column space
$\mtx{A}$ to within precision $\varepsilon$.
Furthermore, we say that an $m\times \ell$ matrix $\mtx{Q}$ is an
$\varepsilon$-\textit{basis matrix} for $\mtx{A}$
if its columns are orthonormal, and if $||\mtx{A} - \mtx{Q}\mtx{Q}^{*}\mtx{A}|| \leq \varepsilon$.
\end{definition}

\subsection{Functions for low-rank factorizations}
For future reference, we introduce two functions ``\texttt{qr}'' and ``\texttt{svd}'' that
can operate in three difference modes. In the first mode, they produce the full
(``economy size'') factorizations described in section \ref{sec:QR} and \ref{sec:SVD}, respectively,
$$
[\mtx{Q},\mtx{R},J] = \texttt{qr}(\mtx{A}),
\qquad
[\mtx{U},\mtx{D},\vct{V}] = \texttt{svd}(\mtx{A}),
\qquad\mbox{and}\qquad
[\mtx{X},J] = \texttt{id}(\mtx{A}).
$$
In practice, we execute these factorizations using standard LAPACK library functions.
In the second mode, we provide an integer $k$ and obtain partial factorizations of
rank $k$,
$$
[\mtx{Q},\mtx{R},J] = \texttt{qr}(\mtx{A},k),
\qquad
[\mtx{U},\mtx{D},\vct{V}] = \texttt{svd}(\mtx{A},k),
\qquad\mbox{and}\qquad
[\mtx{X},J] = \texttt{id}(\mtx{A},k).
$$
Then the matrices $\mtx{Q},\,\mtx{U},\,\mtx{D},\,\mtx{V}$ have precisely $k$
columns, and $\mtx{R}$ has precisely $k$ rows.
In the third mode, we provide a real number $\varepsilon \in (0,1)$ and obtain
partial factorizations
$$
[\mtx{Q},\mtx{R},J] = \texttt{qr}(\mtx{A},\varepsilon),
\qquad
[\mtx{U},\mtx{D},\vct{V}] = \texttt{svd}(\mtx{A},\varepsilon),
\qquad\mbox{and}\qquad
[\mtx{X},J] = \texttt{id}(\mtx{A},\varepsilon),
$$
such that
$$
\|\mtx{A}(\colon,J) - \mtx{Q}\mtx{R}\| \leq \varepsilon
\qquad
\|\mtx{A} - \mtx{U}\mtx{D}\mtx{V}^{*}\| \leq \varepsilon,
\qquad\mbox{and}\qquad
\|\mtx{A} - \mtx{A}(:,J)\mtx{X}\| \leq \varepsilon.
$$
In practice, for a \textit{small} input matrix $\mtx{A}$,
we execute mode 2 and mode 3 by calling the LAPACK routine for a
\textit{full} factorization (the ID can be obtained from the full QR),
and then simply truncating the result. If $\mtx{A}$ is large, then we
use the randomized sampling technique
of Section \ref{sec:randcomp}.

\begin{remark}
The differentiation between modes 2 and 3 for \texttt{qr} and \texttt{svd} is
communicated by whether the second argument is an integer (mode 2) or a real
number $\varepsilon \in (0,1)$ (mode 3). This is slightly questionable notation,
but it keeps the formulas clean, and hopefully does not cause confusion.
\end{remark}

\subsection{A binary tree structure}
\label{sec:tree}

\begin{figure}
\setlength{\unitlength}{1mm}
\fbox{
\begin{picture}(160,41)
\put(14, 0){\includegraphics[height=35mm]{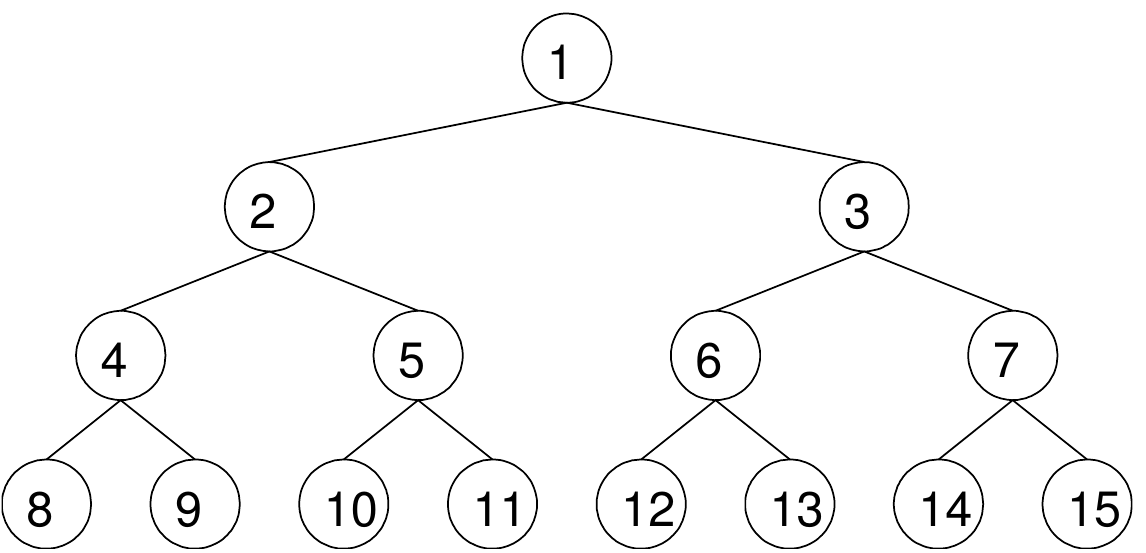}}
\put( 0,31){Level $0$}
\put( 0,21){Level $1$}
\put( 0,11){Level $2$}
\put( 0, 1){Level $3$}
\put(88,31){\footnotesize$I_{1} = [1,\,2,\,\dots,\,400]$}
\put(88,21){\footnotesize$I_{2} = [1,\,2,\,\dots,\,200]$, $I_{3} = [201,\,202,\,\dots,\,400]$}
\put(88,11){\footnotesize$I_{4} = [1,\,2,\,\dots,\,100]$, $I_{5} = [101,\,102,\,\dots,\,200]$, \dots}
\put(88, 1){\footnotesize$I_{8} = [1,\,2,\,\dots,\,50]$, $I_{9} = [51,\,52,\,\dots,\,100]$, \dots}
\end{picture}}
\caption{Numbering of nodes in a fully populated binary tree with $L=3$ levels.
The root is the original index vector $I = I_{1} = [1,\,2,\,\dots,\,400]$.}
\label{fig:tree}
\end{figure}
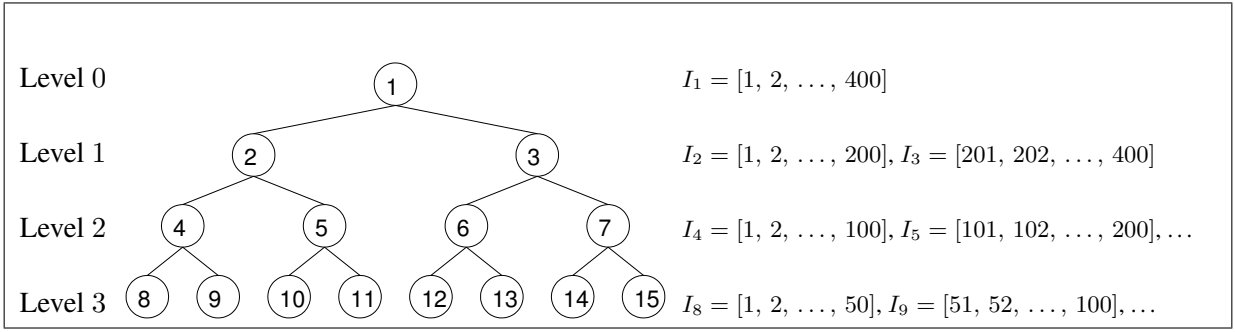

Both the HODLR and the HBS representations
of an $M\times M$ matrix $\mtx{A}$ are
based on a partition of the index vector $I = [1,\,2,\,\dots,\,M]$
into a binary tree structure.
We let $I$ form the root of the tree, and give it the index $1$,
$I_{1} = I$. We next split the root into two roughly equi-sized
vectors $I_{2}$ and $I_{3}$ so that $I_{1} = I_{2} \cup I_{3}$.
The full tree is then formed by continuing to subdivide any interval
that holds more than some preset fixed number $m$ of indices.
We use the integers $\ell = 0,\,1,\,\dots,\,L$ to label the different
levels, with $0$ denoting the coarsest level.
A \textit{leaf} is a node corresponding to a vector that never got split.
For a non-leaf node $\tau$, its \textit{children} are the two boxes
$\sigma_{1}$ and $\sigma_{2}$ such that $I_{\tau} = I_{\sigma_{1}} \cup I_{\sigma_{2}}$,
and $\tau$ is then the \textit{parent} of $\sigma_{1}$ and $\sigma_{2}$.
Two boxes with the same parent are called \textit{siblings}. These
definitions are illustrated in Figure \ref{fig:tree}. For any node $\tau$,
let $n_{\tau}$ denote the number of indices in $I_{\tau}$.

\subsection{The HODLR data sparse matrix format}
\label{sec:introHODLR}
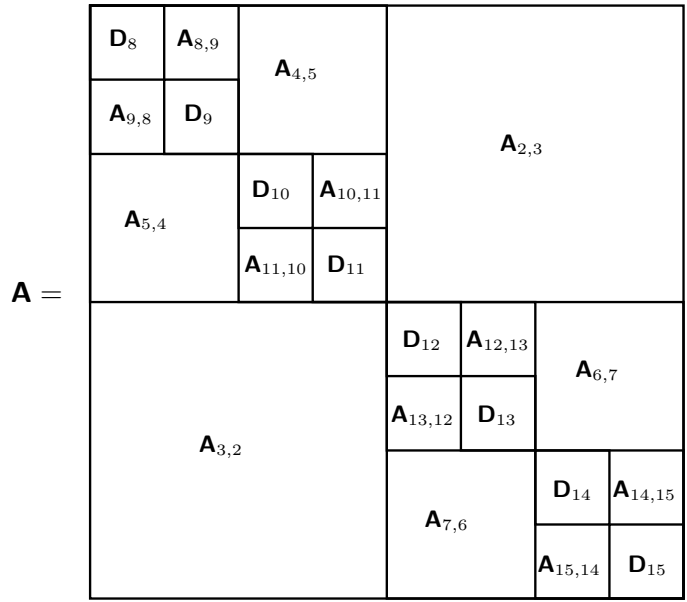
\begin{figure}
\begin{center}
\setlength{\unitlength}{1mm}
\begin{picture}(80,80)
\put(-10,40){$\mtx{A} = $}
\put(00,00){\includegraphics[width=80mm]{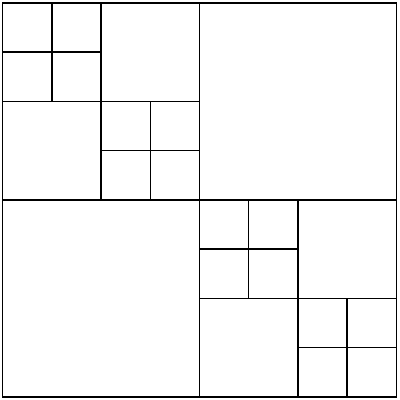}}
\put(55,60){\footnotesize$\mtx{A}_{2,3}$}
\put(15,20){\footnotesize$\mtx{A}_{3,2}$}
\put(25,70){\footnotesize$\mtx{A}_{4,5}$}
\put(05,50){\footnotesize$\mtx{A}_{5,4}$}
\put(65,30){\footnotesize$\mtx{A}_{6,7}$}
\put(45,10){\footnotesize$\mtx{A}_{7,6}$}
\put(03,74){\footnotesize$\mtx{D}_{8}$}
\put(12,74){\footnotesize$\mtx{A}_{8,9}$}
\put(03,64){\footnotesize$\mtx{A}_{9,8}$}
\put(13,64){\footnotesize$\mtx{D}_{9}$}
\put(22,54){\footnotesize$\mtx{D}_{10}$}
\put(31,54){\footnotesize$\mtx{A}_{10,11}$}
\put(21,44){\footnotesize$\mtx{A}_{11,10}$}
\put(32,44){\footnotesize$\mtx{D}_{11}$}
\put(42,34){\footnotesize$\mtx{D}_{12}$}
\put(50.5,34){\footnotesize$\mtx{A}_{12,13}$}
\put(40.5,24){\footnotesize$\mtx{A}_{13,12}$}
\put(52,24){\footnotesize$\mtx{D}_{13}$}
\put(62,14){\footnotesize$\mtx{D}_{14}$}
\put(70,14){\footnotesize$\mtx{A}_{14,15}$}
\put(60,04){\footnotesize$\mtx{A}_{15,14}$}
\put(72,04){\footnotesize$\mtx{D}_{15}$}
\end{picture}
\end{center}
\caption{A HODLR matrix $\mtx{A}$ tesselated in accordance with the tree in Figure \ref{fig:tree}.
Every off-diagonal block $\mtx{A}_{\alpha,\beta}$ that is marked in the figure should
have $\varepsilon$-rank at most $k$.}
\label{fig:HODLR_tessellation}
\end{figure}

The \textit{Hierarchically Off-Diagonal Low Rank (HODLR)} property is, as the
name implies, a condition that the off-diagonal blocks of a matrix $\mtx{A}$
should have low (numerical) rank. To be precise, given a hierarchical partitioning
of the index vector, cf.~Section \ref{sec:tree}, a computational tolerance
$\varepsilon$, and a bound on the rank $k$, we require that for any
sibling pair $\{\alpha,\beta\}$, the corresponding block
$$
\mtx{A}_{\alpha,\beta} = \mtx{A}(I_{\alpha},I_{\beta})
$$
should have $\varepsilon$-rank at most $k$. The tessellation resulting from the
tree in Figure \ref{fig:tree} is shown in Figure \ref{fig:HODLR_tessellation}.
We then represent each off-diagonal block via a rank-$k$ factorization
$$
\begin{array}{cccccccccccccccc}
\mtx{A}_{\alpha,\beta} &=& \mtb{U}_{\alpha} & \tilde{\mtx{A}}_{\alpha,\beta} & \mtb{U}_{\beta},\\
n_{\alpha} \times n_{\beta} && n_{\alpha} \times k & k \times k & k\times n_{\beta}
\end{array}
$$
where $\mtb{U}_{\alpha}$ and $\mtb{B}_{\beta}$ are orthonormal matrices. It is
easily verified that if we required each leaf node to have at most $O(k)$ nodes,
then it takes $O(k\,N\,\log N)$ storage to store all factors required to represent
$\mtx{A}$, and a matrix-vector multiplication can be executed using $O(k\,N\,\log N)$
flops.

For future reference, we define for a given HODLR matrix $\mtx{A}$ a
``level-truncated'' matrix $\mtx{A}^{(\ell)}$ as the matrix obtained
by zeroing out any block associated with levels finer than $\ell$.
In other words,
$$
\mtx{A}^{(1)} =
\begin{array}{|l|l|} \hline
\mtx{0} & \mtx{A}_{2,3} \\ \hline
\mtx{A}_{3,2} & \mtx{0} \\ \hline
\end{array},
\qquad
\mtx{A}^{(2)} =
\begin{array}{|l|l|l|l|}
\hline
\mtx{0}       & \mtx{A}_{4,5} & \multicolumn{2}{c|}{\multirow{2}{*}{$\mtx{A}_{2,3}$}} \\ \cline{1-2}
\mtx{A}_{5,4} &    \mtx{0}    & \multicolumn{2}{c|}{}                           \\ \hline
\multicolumn{2}{|c|}{\multirow{2}{*}{$\mtx{A}_{3,2}$}} &   \mtx{0}                 & \mtx{A}_{6,7}                  \\ \cline{3-4}
\multicolumn{2}{|c|}{}                  & \mtx{A}_{7,6}                   &     \mtx{0}              \\ \hline
\end{array},
\qquad
\mbox{etc.}
$$

\begin{remark}
In our description of rank-structured matrices in sections \ref{sec:introHODLR} and \ref{sec:introHBS},
we generally assume that the numerical rank is the same number $k$ for every off-diagonal block.
In practice, we typically estimate the $\varepsilon$-rank for any specific off-diagonal block
adaptively to save both storage and flop counts.
\end{remark}

\subsection{The HBS data sparse matrix format}
\label{sec:introHBS}
The HODLR format is simple to describe and to use, but is slightly inefficient in that
it requires the user to store for every node $\tau$, the basis matrices $\mtb{U}_{\tau}$
and $\mtb{V}_{\tau}$, which can be quite long. The \textit{Hierarchically Block Separable (HBS)}
format is designed to overcome this problem by expressing these matrices \textit{hierarchically.}
To be precise, suppose that $\tau$ is a node which children $\{\alpha,\beta\}$, and that we
can find a \textit{short} matrix $\mtx{U}_{\tau}$ such that
\begin{equation}
\label{eq:babytelescope}
\begin{array}{cccccccccccccccccccc}
\mtb{U}_{\tau}
&=&
\left[\begin{array}{cc}
\mtb{U}_{\alpha} & \mtx{0} \\
\mtx{0} & \mtb{U}_{\beta}
\end{array}\right]
&
\mtx{U}_{\tau}.\\
n_{\tau} \times k && n_{\tau}\times 2k & 2k \times k
\end{array}
\end{equation}
The point is that if we have the long basis matrices $\mtb{U}_{\alpha}$ and $\mtb{U}_{\beta}$
available for the children, then all we need to store in order to be able to apply
$\mtb{U}_{\tau}$ is the short matrix $\mtx{U}_{\tau}$. This process can now be continued
recursively. For instance, if $\{\gamma,\delta\}$ are the children of $\alpha$,
and $\{\nu,\mu\}$ are the children of $\beta$, we assume there exist matrices
$\mtx{U}_{\alpha}$ and $\mtx{U}_{\beta}$ such that
$$
\begin{array}{cccccccccccccccccccc}
\mtb{U}_{\alpha}
&=&
\left[\begin{array}{cc}
\mtb{U}_{\gamma} & \mtx{0} \\
\mtx{0} & \mtb{U}_{\delta}
\end{array}\right]
&
\mtx{U}_{\alpha},\\
n_{\alpha} \times k && n_{\alpha}\times 2k & 2k \times k
\end{array}
\qquad\mbox{and}\qquad
\begin{array}{cccccccccccccccccccc}
\mtb{U}_{\beta}
&=&
\left[\begin{array}{cc}
\mtb{U}_{\mu} & \mtx{0} \\
\mtx{0} & \mtb{U}_{\nu}
\end{array}\right]
&
\mtx{U}_{\beta}.\\
n_{\beta} \times k && n_{\beta}\times 2k & 2k \times k
\end{array}
$$
Then $\mtb{U}_{\tau}$ can be expressed as
$$
\begin{array}{cccccccccccccccccccc}
\mtb{U}_{\tau}
&=&
\left[\begin{array}{cccc}
\mtb{U}_{\gamma} & \mtx{0} & \mtx{0} & \mtx{0} \\
\mtx{0} & \mtb{U}_{\delta} & \mtx{0} & \mtx{0} \\
\mtx{0} & \mtx{0} & \mtb{U}_{\mu   } & \mtx{0} \\
\mtx{0} & \mtx{0} & \mtx{0} & \mtb{U}_{\nu}
\end{array}\right]
&
\left[\begin{array}{cc}
\mtx{U}_{\alpha} & \mtx{0} \\
\mtx{0} & \mtx{U}_{\beta}
\end{array}\right]
&
\mtx{U}_{\tau}.\\
n_{\tau} \times k && n_{\tau}\times 4k & 4k \times 2k & 2k \times k
\end{array}
$$
By continuing this process down to the leaves, it becomes clear that
we only need to store the ``long'' basis matrices for a leaf node
(and they are not in fact long for a leaf node!); for every other
node, it is sufficient to store the small matrix $\mtx{U}_{\tau}$.
The process for storing the long basis matrices $\mtb{V}_{\tau}$
via small matrices $\mtx{V}_{\tau}$ of size $2k\times k$ is of course
exactly analogous.

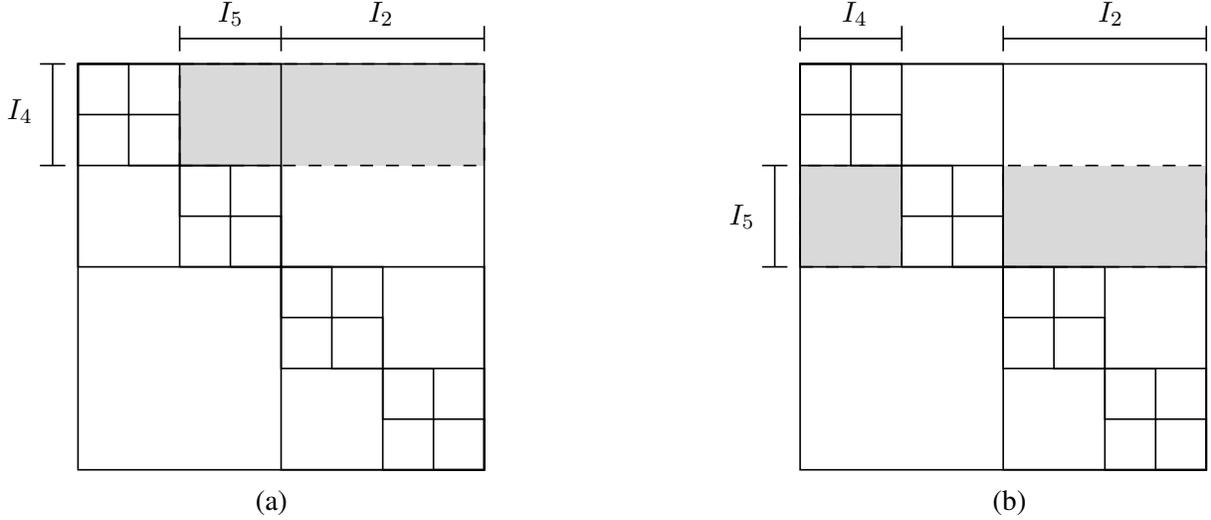
\begin{figure}
\begin{center}
\setlength{\unitlength}{1mm}
\begin{picture}(160,69)
\put(004,05){\includegraphics[height=60mm]{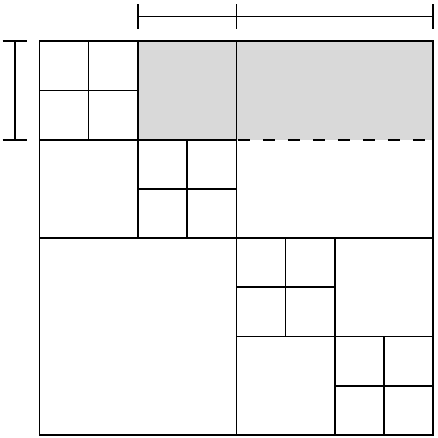}}
\put(000,52){$I_{4}$}
\put(028,65){$I_{5}$}
\put(048,65){$I_{2}$}
\put(033,00){(a)}
\put(100,05){\includegraphics[height=60mm]{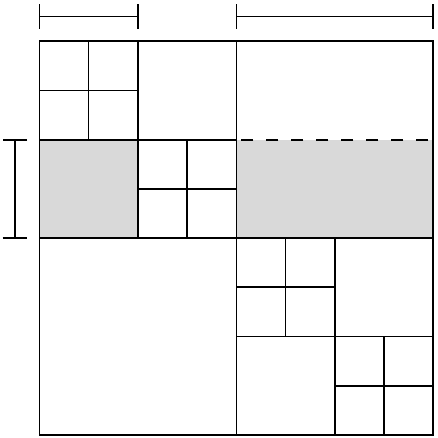}}
\put(096,38){$I_{5}$}
\put(111,65){$I_{4}$}
\put(145,65){$I_{2}$}
\put(131,00){(b)}
\end{picture}
\end{center}
\caption{Illustration of the \textit{neutered row blocks} for the nodes $4$ and $5$, with parent $2$.
(a) The block $\mtx{A}(I_{4},I_{4}^{\rm c})$ is marked in grey.
Observe that $\mtx{A}(I_{4},I_{4}^{\rm c}) = [\mtx{A}(I_{4},I_{5}),\,\mtx{A}(I_{4},I_{2})]$.
(b) The block $\mtx{A}(I_{5},I_{5}^{\rm c})$ is marked in grey.
Observe that $\mtx{A}(I_{5},I_{5}^{\rm c}) = [\mtx{A}(I_{5},I_{4}),\,\mtx{A}(I_{5},I_{2})]$.}
\label{fig:neuteredblocks}
\end{figure}

For a general HODLR matrix, there is no guarantee that a relationship such as
(\ref{eq:babytelescope}) need hold. We need to impose an additional condition
on the long basis matrices $\mtb{U}_{\tau}$ and $\mtb{V}_{\tau}$. To this end,
given a node $\tau$, let us define a \textit{neutered row block} as the off-diagonal
block $\mtx{A}(I_{\tau},I_{\tau}^{\rm c})$, where $I_{\tau}^{\rm c}$ is the
complement of $I_{\tau}$ within the vector $[1,2,3,\dots,N]$, cf.~Figure \ref{fig:neuteredblocks}.
We then require that the columns of the long basis matrix $\mtb{U}_{\tau}$ must span
the columns of $\mtx{A}(I_{\tau},I_{\tau}^{\rm c})$.
Observe that for a node $\tau$ with sibling $\sigma$, the sibling
matrix $\mtx{A}(I_{\tau},I_{\sigma})$ is a submatrix of the neutered row block
$\mtx{A}(I_{\tau},I_{\tau}^{\rm c})$ since $I_{\sigma} \subseteq I_{\tau}^{\rm c}$.
This means that the new requirement of the basis matrices is more restrictive, and
that typically the ranks required will be larger for any given precision. However,
once the long basis matrices satisfy the more restrictive requirement, it is
\textit{necessarily} the case that (\ref{eq:babytelescope}) holds for some
small matrix $\mtx{U}_{\tau}$.

We analogously define
the \textit{neutered column block} for $\tau$ as the matrix
$\mtx{A}(I_{\tau}^{\rm c},I_{\tau})$ and require that the columns of
$\mtb{V}_{\tau}$ span the rows of $\mtx{A}(I_{\tau}^{\rm c},I_{\tau})$.

\begin{definition}
We say that a HODLR matrix $\mtx{A}$ is an HBS matrix if,
for every parent node $\tau$ with children $\{\alpha,\beta\}$,
there exist ``small'' basis matrices $\mtx{U}_{\tau}$ and $\mtx{V}_{\tau}$
such that
$$
\begin{array}{cccccccccccccccccccc}
\mtb{U}_{\tau}
&=&
\left[\begin{array}{cc}
\mtb{U}_{\alpha} & \mtx{0} \\
\mtx{0} & \mtb{U}_{\beta}
\end{array}\right]
&
\mtx{U}_{\tau},\\
n_{\tau} \times k && n_{\tau}\times 2k & 2k \times k
\end{array}
\qquad\mbox{and}\qquad
\begin{array}{cccccccccccccccccccc}
\mtb{V}_{\tau}
&=&
\left[\begin{array}{cc}
\mtb{V}_{\alpha} & \mtx{0} \\
\mtx{0} & \mtb{V}_{\beta}
\end{array}\right]
&
\mtx{V}_{\tau}.\\
n_{\tau} \times k && n_{\tau}\times 2k & 2k \times k
\end{array}
$$
\end{definition}

While standard practice is to require all basis matrices $\mtb{U}_{\tau}$ and $\mtb{V}_{\tau}$
to be \textit{orthonormal}, we have found that it is highly convenient to use the
interpolatory decomposition (ID) to represent the off-diagonal blocks. The key advantage
is that then the sibling interaction matrices which be \textit{submatrices} of the original
matrix. This improves interpretability, and also slightly reduces storage requirements.

\begin{definition}
\label{def:HBSID}
We say that a HBS matrix $\mtx{A}$ is in HBS-ID format if every
basis matrix $\mtx{U}_{\tau}$ and $\mtx{V}_{\tau}$ contains
a $k\times k$ identity matrix, and every siblin interaction
matrix $\tilde{\mtx{A}}_{\alpha,\beta}$ is a submatrix of $\mtx{A}$.
In other words, there exist some index sets $\tilde{I}_{\alpha}$
and $\hat{I}_{\beta}$ such that
$$
\tilde{\mtx{A}}_{\alpha,\beta} = \mtx{A}(\tilde{I}_{\alpha},\hat{I}_{\beta}).
$$
The index sets $\tilde{I}_{\tau}$ and $\hat{I}_{\tau}$ are called
the \textit{row skeleton} and \textit{column skeleton} of box $\tau$, respectively.
We enforce that these are ``nested,'' which is to say that if the
children of $\tau$ are $\{\alpha,\beta\}$, then
$$
\tilde{I}_{\tau} \subseteq \tilde{I}_{\alpha} \cup \tilde{I}_{\beta}
\qquad\mbox{and}\qquad
\hat{I}_{\tau} \subseteq \hat{I}_{\alpha} \cup \hat{I}_{\beta}.
$$
\end{definition}

\begin{remark}
\label{remark:spanningskel}
The straight-forward way to build a basis matrix $\mtb{U}_{\tau}$ for a node
$\tau$ is to explicitly form the corresponding neutered row block
$\mtx{A}(I_{\tau},I_{\tau}^{\rm c})$ and then compress it (perform column
pivtoed Gram-Schmidt on its \textit{columns} to form an ON basis $\mtb{U}_{\tau}$,
or perform column pivoted Gram-Schmidt on its \textit{rows} to form the ID).
However, suppose that we can somehow construct a smaller matrix $\mtx{Y}_{\tau}$
of size $n_{\tau} \times \ell$ with the property that the columns of $\mtx{Y}_{\tau}$
span the columns of $\mtx{A}(I_{\tau},I_{\tau}^{\rm c})$. Then it would be
sufficient to process the columns of $\mtx{Y}_{\tau}$ to build a basis for
$\mtx{A}(I_{\tau},I_{\tau}^{\rm c})$. For instance, if we orthonormalize the
columns of $\mtx{Y}_{\tau}$ to form a basis matrix $\mtb{U}_{\tau} = \texttt{qr}(\mtx{Y}_{\tau})$,
then the columns of $\mtb{U}_{\tau}$ will necessarily form an orthonormal basis
for the columns of $\mtx{A}(I_{\tau},I_{\tau}^{\rm c})$. The key point here
is that one can often find such a matrix $\mtx{Y}_{\tau}$ with a small number
$\ell$ of columns. In \cite{2005_martinsson_fastdirect} we use a representation
theorem from potential theory to find such a matrix $\mtx{Y}_{\tau}$ when
$\mtx{A}$ comes from the discretization of a boundary integral equation of
mathematical physics. In \cite{2011_martinsson_randomhudson} we do this
via randomized sampling, so that $\mtx{Y}_{\tau} = \mtx{A}(I_{\tau},I_{\tau}^{\rm c})\,\mtx{\Omega}$
for some Gaussian random matrix $\mtx{\Omega}$. The main point of \cite{2011_martinsson_randomhudson}
is that this can be done by applying all of $\mtx{A}$ to a \textit{single} random matrix with $N\times \ell$
columns, where $\ell \approx k$. In the current manuscript, we use a similar strategy,
but we now require the application of $\mtx{A}$ to a set of $O(k)$ random vector
for each level.
\end{remark}


\section{An algorithm for compressing a HODLR matrix}
\label{sec:HODLR}

The algorithm consists of a single sweep through the levels in the tree of nodes,
starting from the root (the entire domain), and processing each level of successively
smaller blocks at a time.

In presenting the algorithm, we let $\mtx{G}_{\tau}$ denote a random matrix
of size $n_{\tau} \times r$ drawn from a Gaussian distribution. Blocked matrices
are drawn with the blocks to be processed set in red type. To minimize clutter,
blocks of the matrix that will not play a part in the current step are marked
with a star (``*'') and may or may not be zero.

\vspace{2mm}

\noindent
\textit{Processing level 0 (the root of the tree):} Let $\{\alpha,\beta\}$ denote the children
of the root node (in our standard ordering, $\alpha=2$ and $\beta=3$).
Our objective is now to find low-rank factorizations of the off-diagonal blocks
$$
\mtx{A} =
\left[\begin{array}{cc}
* & {\color{red}\mtx{A}_{\alpha\beta}}\\
{\color{red}\mtx{A}_{\alpha\beta}} & *\\
\end{array}\right].
$$
To this end, we build two random matrix, each of size $N\times r$, and defined by
$$
\mtx{\Omega}_{1} = \vtwo{\mtx{G}_{\alpha}}{\mtx{0}}
\quad\mbox{and}\qquad
\mtx{\Omega}_{2} = \vtwo{\mtx{0}}{\mtx{G}_{\beta }}.
$$
Then construct the matrices of samples
$$
\mtx{Y}_{1} = \mtx{A}  \mtx{\Omega_{2}} = \vtwo{\mtx{A}_{\alpha\beta}\mtx{G}_{\beta }}{*}
\qquad\mbox{and}\qquad
\mtx{Y}_{2} = \mtx{A}  \mtx{\Omega_{1}} = \vtwo{*}{\mtx{A}_{\beta\alpha}\mtx{G}_{\alpha}}.
$$
Supported by the results on randomized sampling described in Section \ref{sec:randcomp}, we now know that
it is almost certain that the columns in the top block of $\mtx{Y}_{1}$ will span the columns
of $\mtx{A}_{\alpha\beta}$. By orthonormalizing the columns of $\mtx{Y}_{1}(I_{\alpha},:)$,
we therefore obtain an ON-basis for the column space of $\mtx{A}_{\alpha\beta}$. In other words,
we set
$$
\mtb{U}_{\alpha} = \texttt{qr}(\mtx{Y}_{1}(I_{\alpha},:)),
\quad\mbox{and}\qquad
\mtb{U}_{\beta} = \texttt{qr}(\mtx{Y}_{2}(I_{\beta},:)),
$$
and then we know that $\mtb{U}_{\alpha}$ and $\mtb{U}_{\beta}$ will serve as the relevant
basis matrices in the HODLR representation of $\mtx{A}$. To compute
$\mtb{V}_{\alpha}$, $\mtb{V}_{\beta}$, $\mtx{B}_{\alpha\beta}$ and $\mtx{B}_{\beta\alpha}$
we now need to form the matrices $\mtb{U}_{\alpha}^{*}\mtx{A}_{\alpha\beta}$ and
$\mtb{U}_{\beta}^{*}\mtx{A}_{\beta\alpha}$. To do this through our ``black-box''
matrix-matrix multiplier, we form the new test matrices
$$
\mtx{\Omega}_{1} = \vtwo{\mtb{U}_{\alpha}}{\mtx{0}}
\quad\mbox{and}\qquad
\mtx{\Omega}_{2} = \vtwo{\mtx{0}}{\mtb{U}_{\beta}}.
$$
Then construct the matrices of samples
$$
\mtx{Z}_{1} = \mtx{A}^{*} \mtx{\Omega_{2}} = \vtwo{\mtx{A}_{\beta\alpha}^{*}\mtb{U}_{\beta }}{*}
\qquad\mbox{and}\qquad
\mtx{Z}_{2} = \mtx{A}^{*} \mtx{\Omega_{1}} = \vtwo{*}{\mtx{A}_{\alpha\beta}^{*}\mtb{U}_{\alpha}}.
$$
All that remains is now to compute QR factorizations
\begin{equation}
\label{eq:QRtogetV}
[\mtb{V}_{\alpha},\mtx{B}_{\beta\alpha}] = \texttt{qr}(\mtx{Z}_{1}(I_{\alpha},:)),
\quad\mbox{and}\qquad
[\mtb{V}_{\beta },\mtx{B}_{\alpha\beta}] = \texttt{qr}(\mtx{Z}_{2}(I_{\beta},:)).
\end{equation}

\begin{remark}
At a slight increase in cost, one can obtain \textit{diagonal} sibling
interaction matrices $\mtx{B}_{\alpha\beta}$ and $\mtx{B}_{\beta\alpha}$.
We would then replace the QR factorization in (\ref{eq:QRtogetV}) by
a full SVD
\begin{equation}
\label{eq:SVDtogetV}
[\mtb{V}_{\alpha},\mtx{B}_{\beta\alpha},\hat{\mtx{U}}_{\beta}] = \texttt{svd}(\mtx{Z}_{1}(I_{\alpha},:)),
\quad\mbox{and}\qquad
[\mtb{V}_{\beta },\mtx{B}_{\alpha\beta},\hat{\mtx{U}}_{\alpha}] = \texttt{svd}(\mtx{Z}_{2}(I_{\beta},:)).
\end{equation}
This step then requires an update to the long basis matrices for the column space:
\begin{equation}
\label{eq:SVDtogetB}
\mtb{U}_{\alpha} \leftarrow \mtb{U}_{\alpha}\hat{\mtx{U}}_{\alpha}
\quad\mbox{and}\qquad
\mtb{U}_{\beta } \leftarrow \mtb{U}_{\beta }\hat{\mtx{U}}_{\beta }.
\end{equation}
\end{remark}

\lsp

\noindent
\textit{Processing level 1:}
Now that all off-diagonal blocks on level 1 have been computed, we can use
this information to compress the blocks on level 2. We let $\{\alpha,\beta,\gamma,\delta\}$
denote the boxes on level 2 (in standard ordering, $\alpha=4$, $\beta=5$, $\gamma=6$, $\delta=7$).
Our objective is now to construct low-rank approximations to the blocks marked in red:
$$
\mtx{A} =
\begin{array}{|l|l|l|l|}
\hline
*                  & {\color{red}\mtx{A}_{\alpha\beta}} & \multicolumn{2}{c|}{\multirow{2}{*}{*}} \\ \cline{1-2}
{\color{red}\mtx{A}_{\beta\alpha}}                   &    *                & \multicolumn{2}{c|}{}                  \\ \hline
\multicolumn{2}{|c|}{\multirow{2}{*}{*}} &   *                 & {\color{red}\mtx{A}_{\gamma\delta}}                  \\ \cline{3-4}
\multicolumn{2}{|c|}{}                  & {\color{red}\mtx{A}_{\delta\gamma}}                   &     *              \\ \hline
\end{array}
$$
First, observe that
$$
\mtx{A} - \mtx{A}^{(1)} =
\begin{array}{|l|l|l|l|}
\hline
*                  & {\color{red}\mtx{A}_{\alpha\beta}} & \multicolumn{2}{c|}{\multirow{2}{*}{$\mtx{0}$}} \\ \cline{1-2}
{\color{red}\mtx{A}_{\beta\alpha}}                   &    *                & \multicolumn{2}{c|}{}                  \\ \hline
\multicolumn{2}{|c|}{\multirow{2}{*}{$\mtx{0}$}} &   *                 & {\color{red}\mtx{A}_{\gamma\delta}}                  \\ \cline{3-4}
\multicolumn{2}{|c|}{}                  & {\color{red}\mtx{A}_{\delta\gamma}}                   &     *              \\ \hline
\end{array}
$$
We then define two random test matrices, each of size $N\times r$, via
$$
\mtx{\Omega}_{1} = \vfour{\mtx{G}_{\alpha}}{\mtx{0}}{\mtx{G}_{\gamma}}{\mtx{0}}
\qquad\mbox{and}\qquad
\mtx{\Omega}_{2} = \vfour{\mtx{0}}{\mtx{G}_{\beta}}{\mtx{0}}{\mtx{G}_{\delta}}
$$
We compute the sample matrices via
\begin{equation}
\label{eq:gretchen}
\mtx{Y}_{1} =
\mtx{A}\mtx{\Omega_{2}} - \mtx{A}^{(1)}\mtx{\Omega_{2}} =
\vfour{\mtx{A}_{\alpha\beta}\mtx{G}_{\beta}}{*}{\mtx{A}_{\gamma\delta}\mtx{G}_{\delta}}{*}
\qquad\mbox{and}\qquad
\mtx{Y}_{2} =
\mtx{A}\mtx{\Omega_{1}} - \mtx{A}^{(1)}\mtx{\Omega_{1}} =
\vfour{*}{\mtx{A}_{\beta\alpha}\mtx{G}_{\alpha}}{*}{\mtx{A}_{\delta\gamma}\mtx{G}_{\gamma}}.
\end{equation}
In evaluating $\mtx{Y}_{1}$ and $\mtx{Y}_{2}$, we use the black-box multiplier to
form $\mtx{A}\mtx{\Omega_{2}}$ and $\mtx{A}\mtx{\Omega_{1}}$, and the compressed
representation of $\mtx{A}^{(1)}$ obtained on the previous level to evaluate
$\mtx{A}^{(1)}\mtx{\Omega_{2}}$ and $\mtx{A}^{(1)}\mtx{\Omega_{1}}$.
We now obtain orthonormal bases for the column spaces of the four sibling interaction
matrices by orthonormalizing the pertinent blocks of $\mtx{Y}_{1}$ and $\mtx{Y}_{2}$:
$$
\mtb{U}_{\alpha} = \texttt{qr}(\mtx{Y}_{1}(I_{\alpha},:)),
\quad
\mtb{U}_{\beta } = \texttt{qr}(\mtx{Y}_{2}(I_{\beta },:)),
\quad
\mtb{U}_{\gamma} = \texttt{qr}(\mtx{Y}_{1}(I_{\gamma},:)),
\quad
\mtb{U}_{\delta} = \texttt{qr}(\mtx{Y}_{2}(I_{\delta},:)).
$$
It remains to construct the ON-bases for the corresponding row spaces and the
compressed sibling interaction matrices. To this end, we form two new test matrices,
both of size $N \times r$, via
$$
\mtx{\Omega}_{1} = \vfour{\mtb{U}_{\alpha}}{\mtx{0}}{\mtb{U}_{\gamma}}{\mtx{0}}
\qquad\mbox{and}\qquad
\mtx{\Omega}_{2} = \vfour{\mtx{0}}{\mtb{U}_{\beta}}{\mtx{0}}{\mtb{U}_{\delta}}.
$$
Then the sample matrices are computed via
$$
\mtx{Z}_{1} =
\mtx{A}^{*}\mtx{\Omega_{2}} - \bigl(\mtx{A}^{(1)}\bigr)^{*}\mtx{\Omega_{2}} =
\vfour{\mtx{A}_{\beta \alpha}^{*}\mtb{U}_{\beta }}{*}{\mtx{A}_{\delta\gamma}^{*}\mtb{U}_{\delta}}{*}
\qquad\mbox{and}\qquad
\mtx{Z}_{2} =
\mtx{A}^{*}\mtx{\Omega_{1}} - \bigl(\mtx{A}^{(1)}\bigr)^{*}\mtx{\Omega_{1}} =
\vfour{*}{\mtx{A}_{\alpha\beta }^{*}\mtb{U}_{\alpha}}{*}{\mtx{A}_{\gamma\delta}^{*}\mtb{U}_{\gamma}}.
$$
We obtain diagonal compressed sibling interaction matrices by taking a sequence of dense SVDs of the
relevant sub-blocks, cf.~(\ref{eq:SVDtogetV}),
\begin{align*}
[\mtb{V}_{\alpha},\mtx{B}_{\beta \alpha},\hat{\mtx{U}}_{\beta }] =&\ \texttt{svd}(\mtx{Z}_{1}(I_{\alpha},:)),\\
[\mtb{V}_{\beta },\mtx{B}_{\alpha\beta },\hat{\mtx{U}}_{\alpha}] =&\ \texttt{svd}(\mtx{Z}_{2}(I_{\beta },:)),\\
[\mtb{V}_{\gamma},\mtx{B}_{\delta\gamma},\hat{\mtx{U}}_{\delta}] =&\ \texttt{svd}(\mtx{Z}_{1}(I_{\gamma},:)),\\
[\mtb{V}_{\delta},\mtx{B}_{\gamma\delta},\hat{\mtx{U}}_{\gamma}] =&\ \texttt{svd}(\mtx{Z}_{2}(I_{\delta},:)).
\end{align*}
Finally update the bases for the column-spaces, cf (\ref{eq:SVDtogetB}),
$$
\mtb{U}_{\alpha} \leftarrow \mtb{U}_{\alpha}\hat{\mtx{U}}_{\alpha},
\qquad
\mtb{U}_{\beta } \leftarrow \mtb{U}_{\beta }\hat{\mtx{U}}_{\beta },
\qquad
\mtb{U}_{\gamma} \leftarrow \mtb{U}_{\gamma}\hat{\mtx{U}}_{\gamma},
\qquad
\mtb{U}_{\delta} \leftarrow \mtb{U}_{\delta}\hat{\mtx{U}}_{\delta}.
$$

\lsp

\noindent
\textit{Processing levels 2 through $L-1$:} The processing of every level proceeds in a manner
completely analogous to the processing of level $1$. The relevant formulas are given in Figure
\ref{alg:HODLR}.

\lsp

\noindent
\textit{Processing the leaves:} Once all $L$ levels have been traversed using the procedure
described, compressed representations of all off-diagonal blocks will have been computed. At
this point, all that remains is to extract the diagonal blocks. We illustrate the process
for a simplistic example of a tree with only $L=2$ levels (beyond the root). Letting
$\{\alpha,\beta,\gamma,\delta\}$ denote the leaf nodes, our task is then to extract the
blocks marked in red:
$$
\mtx{A} =
\begin{array}{|l|l|l|l|}
\hline
{\color{red}\mtx{D}_{\alpha}} & * & \multicolumn{2}{c|}{\multirow{2}{*}{*}}  \\ \cline{1-2}
* & {\color{red}\mtx{D}_{\beta }} & \multicolumn{2}{c|}{}                    \\ \hline
\multicolumn{2}{|c|}{\multirow{2}{*}{*}} & {\color{red}\mtx{D}_{\gamma}} & * \\ \cline{3-4}
\multicolumn{2}{|c|}{}                   & * & {\color{red}\mtx{D}_{\delta}} \\ \hline
\end{array}
=
\begin{array}{|l|l|l|l|}
\hline
{\color{red}\mtx{D}_{\alpha}} & \mtx{0} & \multicolumn{2}{c|}{\multirow{2}{*}{\mtx{0}}}  \\ \cline{1-2}
\mtx{0} & {\color{red}\mtx{D}_{\beta }} & \multicolumn{2}{c|}{}                    \\ \hline
\multicolumn{2}{|c|}{\multirow{2}{*}{\mtx{0}}} & {\color{red}\mtx{D}_{\gamma}} & \mtx{0} \\ \cline{3-4}
\multicolumn{2}{|c|}{}                   & \mtx{0} & {\color{red}\mtx{D}_{\delta}} \\ \hline
\end{array} + \mtx{A}^{(2)}.
$$
Since the diagonal blocks are not rank-deficient, we will extract them directly, without using
randomized sampling. To describe the process, we assume at first (for simplicity) that every
diagonal block has the same size, $m\times m$. We then choose a test matrix of size $N\times m$
$$
\mtx{\Omega} = \vfour{\mtx{I}_{m}}{\mtx{I}_{m}}{\mtx{I}_{m}}{\mtx{I}_{m}},
$$
and trivially extract the diagonal blocks via the sampling
$$
\mtx{Y} =
\mtx{A}\mtx{\Omega} - \mtx{A}^{(2)}\mtx{\Omega} =
\vfour{\mtx{D}_{\alpha}}{\mtx{D}_{\beta }}{\mtx{D}_{\gamma}}{\mtx{D}_{\delta}}.
$$
The diagonal blocks can then be read off directly from $\mtx{Y}$.

For the general case where the leaves may be of different sizes, we form a test matrix
$\mtx{\Omega}$ of size $N\times m$, where $m = \max\{n_{\tau}\,\colon\,\tau\mbox{ is a leaf}\}$,
such that
$$
\mtx{\Omega}(I_{\tau},:) = \bigl[\mtx{I}_{n_{\tau}}\ \texttt{zeros}(n_{\tau},m-n_{\tau})\bigr].
$$
In other words, we simply pad a few zero columns at the end.

The entire algorithm is summarized in Figure \ref{alg:HODLR}.

\begin{figure}
\fbox{\begin{minipage}{190mm}
      \begin{tabbing}
      \hspace{5mm} \= \hspace{5mm} \= \hspace{5mm} \= \hspace{5mm}\kill
      \textit{\color{blue}Build compressed representations of all off-diagonal blocks.}\\
      \textbf{loop} over levels $\ell = 0:(L-1)$\\[2mm]
      \> \textit{\color{blue}Build the random matrices $\mtx{\Omega}_{1}$ and $\mtx{\Omega}_{2}$.}\\
      \> $\mtx{\Omega}_{1} = \texttt{zeros}(n,r)$\\
      \> $\mtx{\Omega}_{2} = \texttt{zeros}(n,r)$\\
      \> \textbf{loop} over boxes $\tau$ on level $\ell$\\
      \> \> Let $\{\alpha,\beta\}$ denote the children of box $\tau$.\\
      \> \> $\mtx{\Omega}_{1}(I_{\alpha},:) = \texttt{randn}(n_{\alpha},r)$\\
      \> \> $\mtx{\Omega}_{2}(I_{\beta },:) = \texttt{randn}(n_{\beta },r)$\\
      \> \textbf{end loop}\\[2mm]
      \> \textit{\color{blue}Apply $\mtx{A}$ to build the samples for the incoming basis matrices.}\\
      \> $\mtx{Y}_{1} = \mtx{A}\mtx{\Omega}_{2} - \mtx{A}^{(\ell)}\mtx{\Omega}_{2}$\\
      \> $\mtx{Y}_{2} = \mtx{A}\mtx{\Omega}_{1} - \mtx{A}^{(\ell)}\mtx{\Omega}_{1}$\\[2mm]
      \> \textit{\color{blue}Orthonormalize the sample matrices to build the incoming basis matrices.}\\
      \> \textbf{loop} over boxes $\tau$ on level $\ell$\\
      \> \> Let $\{\alpha,\beta\}$ denote the children of box $\tau$.\\
      \> \> $\mtb{U}_{\alpha} = \texttt{qr}(\mtx{Y}_{1}(I_{\alpha},:))$.\\
      \> \> $\mtb{U}_{\beta } = \texttt{qr}(\mtx{Y}_{2}(I_{\beta },:))$.\\
      \> \> $\mtx{\Omega}_{1}(I_{\alpha},:) = \mtb{U}_{\alpha}$\\
      \> \> $\mtx{\Omega}_{2}(I_{\beta },:) = \mtb{U}_{\beta }$\\
      \> \textbf{end loop}\\[2mm]
      \> \textit{\color{blue}Apply $\mtx{A}^{*}$ to build the samples for the outgoing basis matrices.}\\
      \> $\mtx{Z}_{1} = \mtx{A}^{*}\mtx{\Omega}_{2} - \bigl(\mtx{A}^{(\ell)}\bigr)^{*}\mtx{\Omega}_{2}$\\
      \> $\mtx{Z}_{2} = \mtx{A}^{*}\mtx{\Omega}_{1} - \bigl(\mtx{A}^{(\ell)}\bigr)^{*}\mtx{\Omega}_{1}$\\[2mm]
      \> \textit{\color{blue}Take local SVDs to build incoming basis matrices and sibling interaction matrices.}\\
      \> \textit{\color{blue}We determine the actual rank, and update the $\mtb{U}_{*}$ basis matrices accordingly.}\\
      \> \textbf{loop} over boxes $\tau$ on level $\ell$\\
      \> \> Let $\{\alpha,\beta\}$ denote the children of box $\tau$.\\
      \> \> $[\mtb{V}_{\alpha},\mtx{B}_{\beta\alpha},\hat{\mtx{U}}_{\beta }] = \texttt{svd}(\mtx{Z}_{1}(I_{\alpha},:),\varepsilon)$.\\
      \> \> $[\mtb{V}_{\beta },\mtx{B}_{\alpha\beta},\hat{\mtx{V}}_{\alpha}] = \texttt{svd}(\mtx{Z}_{2}(I_{\beta},:),\varepsilon)$.\\
      \> \> $\mtb{U}_{\beta } \leftarrow \mtb{U}_{\beta }\hat{\mtx{U}}_{\beta }$.\\
      \> \> $\mtb{U}_{\alpha} \leftarrow \mtb{U}_{\alpha}\hat{\mtx{U}}_{\alpha}$.\\
      \> \textbf{end loop}\\
      \textbf{end loop}\\[2mm]
      \textit{\color{blue}Extract the diagonal matrices.}\\
      $n_{\rm max} = \texttt{max}\,\{n_{\tau}\, \colon\, \tau\mbox{ is a leaf}\}$\\
      $\mtx{\Omega} = \texttt{zeros}(N,n_{\rm max})$\\
      \textbf{loop} over leaf boxes $\tau$\\
      \> $\mtx{\Omega}(I_{\tau},1:n_{\tau}) = \texttt{eye}(n_{\tau})$.\\
      \textbf{end loop}\\
      $\mtx{Y} = \mtx{A}\mtx{\Omega} - \mtx{A}^{(L)}\mtx{\Omega}$\\
      \textbf{loop} over leaf boxes $\tau$\\
      \> $\mtx{D}_{\tau} = \mtx{Y}(I_{\tau},1:n_{\tau})$.\\
      \textbf{end loop}
      \end{tabbing}
      \end{minipage}}
\caption{Randomized compression of a HODLR matrix.}
\label{alg:HODLR}
\end{figure}

\subsection{Asymptotic complexity}
\label{sec:HODLRasympt}
Let $L$ denote the number of levels in the tree.
We find that $L\sim \log N$. Let $T_{\rm flop}$ denote the time required for a flop,
let $T_{\rm mult}$ denote the time required to apply $\mtx{A}$ or $\mtx{A}^{*}$ to
a vector, and let $T_{\mtx{A}^{(\ell)}}$ denote the time required to apply $\mtx{A}^{(\ell)}$
to a vector. Then the cost of processing level $\ell$ is
\begin{equation}
\label{eq:dev1}
T_{\ell} \sim T_{\rm mult} \times k + T_{\rm flop} \times 2^{\ell}\,k\,\frac{N}{2^{\ell}} + T_{\mtx{A}^{(\ell)}}\times k,
\end{equation}
since on level $\ell$ there are $2^{\ell}$ blocks to be processed, and each ``long'' matrix at this
level has height $2^{-\ell}\,N$ and width $O(k)$. Further, we find that the cost of
applying $\mtx{A}^{(\ell)}$ to a single vector is
\begin{equation}
\label{eq:dev2}
T_{\mtx{A}^{(\ell)}}
\sim T_{\rm flop} \times \sum_{j=0}^{\ell} 2^{j}\,k\,\frac{N}{2^{j}}
\sim T_{\rm flop} \times \ell\,k\,N.
\end{equation}
Combining (\ref{eq:dev1}) and (\ref{eq:dev2}) and summing from $\ell=0$ to $\ell=L$, we find (using that $L \sim \log N$)
\begin{equation}
\label{eq:T_HODLR}
T_{\rm compress} \sim T_{\rm mult} \times k\,\log N + T_{\rm flop} \times k^{2}\,N\,\bigl(\log N\bigr)^{2}.
\end{equation}

\section{An algorithm for compressing an HBS matrix}

The algorithm for computing a compressed representation of an HBS matrix is a slight
variation of the algorithm for a HODLR matrix described in Section \ref{sec:HODLR}.
The key difference is that
the long basis matrices $\mtb{U}_{\alpha}$ and $\mtb{V}_{\alpha}$ associated with any
node now must satisfy a more rigorous requirement. We will accomplish this objective
by constructing for every node $\alpha$, two new ``long'' sampling matrices
$\mtb{Y}_{\alpha}$ and $\mtb{Z}_{\alpha}$, each of size $n_{\alpha}\times r$ that help
transmit information from the higher levels to the node $\alpha$. The presentation will
start in Section \ref{sec:basicHBS} with a description of the modification to the
scheme of \ref{sec:HODLR} required to enforce the more rigorous requirement. In this
initial description, we will assume that all four ``long'' matrices
associated with a node ($\mtb{U}_{\tau}$, $\mtb{V}_{\tau}$, $\mtb{Y}_{\tau}$, $\mtb{Z}_{\tau}$)
are stored explicitly, resulting in an $O(k\,N\,\log N)$ memory requirement, just like for
the HODLR algorithm. In Section \ref{sec:memoryHBS} we show that while these ``long'' matrices
do need to be temporarily built and processed, they can be stored \textit{implicitly,} which
will allow the algorithm to use only $O(k\,N)$ memory. Finally, Section \ref{sec:HBSID} will
describe how to construct a representation using interpolatory decompositions in all low-rank
approximations.

\subsection{A basic scheme for compressing an HBS matrix}
\label{sec:basicHBS}
Throughout this section, let $\alpha$ denote a node with a parent $\tau$ that is not the root,
and with a sibling $\beta$. We will first describe the process for building the long basis
matrices $\{\mtb{U}_{\tau}\}_{\tau}$. To this end, recall that the difference in requirements
on the long basis matrices in the two frameworks is as follows:
\begin{align*}
\mbox{HODLR framework:}&\ \mbox{The columns of $\mtb{U}_{\alpha}$ need to span the columns of }\ \mtx{A}(I_{\alpha},I_{\beta})\\
\mbox{HBS framework:}  &\ \mbox{The columns of $\mtb{U}_{\alpha}$ need to span the columns of }\ \mtx{A}(I_{\alpha},I_{\alpha}^{\rm c}).
\end{align*}
The assertion that the requirements on a basis in the HBS framework is more stringent follows
from the fact that $I_{\beta} \subseteq I_{\alpha}^{\rm c}$, and that $I_{\beta}$ is in general
much smaller than $I_{\alpha}^{\rm c}$. Now note that
\begin{equation}
\label{eq:blabla}
\mtx{A}(I_{\alpha},I_{\alpha}^{\rm c}) = \bigl[\mtx{A}(I_{\alpha},I_{\beta})\ \mtx{A}(I_{\alpha},I_{\tau}^{\rm c})\bigr]\,\mtx{P}.
\end{equation}
In (\ref{eq:blabla}), the matrix $\mtx{P}$ is a permutation matrix whose effect is to reorder the columns.
For purposes of constructing a basis for the column space, the matrix $\mtx{P}$ can be ignored.
The idea is now to introduce a new sampling matrix $\mtb{Y}_{\alpha}$ of size $n_{\alpha} \times r$ that encodes all the information
that needs to be transmitted from the parent $\tau$. Specifically, we ask that:
$$
\mbox{The columns of $\mtb{Y}_{\alpha}$ span the columns of $\mtx{A}(I_{\alpha},I_{\tau}^{\rm c})$ (to within precision $\varepsilon$).}
$$
Then, when processing box $\alpha$, we will sample $\mtx{A}(I_{\alpha},I_{\beta})$ using a
Gaussian matrix $\mtx{G}_{\beta}$ of size $n_{\beta} \times r$ just as in the HODLR algorithm.
In the end, we will build $\mtb{U}_{\alpha}$ by combining the two sets of samples
$$
[\mtb{U}_{\alpha},\mtx{D}_{\alpha},\sim] =
\texttt{svd}\bigl(\bigl[\mtx{A}(I_{\alpha},I_{\beta})\mtx{G}_{\beta},\ \mtb{Y}_{\alpha}\bigr],r\bigr).
$$
In other words, we take a matrix $\bigl[\mtx{A}(I_{\alpha},I_{\beta})\mtx{G}_{\beta},\ \mtb{Y}_{\alpha}\bigr]$
of size $n_{\alpha}\times 2r$ and extract its leading $r$ singular components (the trailing $r$
components are discarded). All that remains is to build the sample matrices $\mtb{Y}_{\gamma}$ and
$\mtb{Y}_{\delta}$ that transmit information to the children $\{\gamma,\delta\}$ of $\alpha$. To be
precise, let $J_{\gamma}$ and $J_{\delta}$ denote the local \textit{relative} index vectors, so that
$$
I_{\gamma} = I_{\alpha}(J_{\gamma})
\qquad\mbox{and}\qquad
I_{\delta} = I_{\alpha}(J_{\delta}).
$$
Then set
$$
\mtb{Y}_{\gamma} = \mtb{U}(J_{\gamma},:)\mtx{D}_{\alpha}
\qquad\mbox{and}\qquad
\mtb{Y}_{\delta} = \mtb{U}(J_{\delta},:)\mtx{D}_{\alpha}.
$$

The process for building the long basis matrices $\{\mtb{V}_{\alpha}\}_{\alpha}$
is entirely analogous to the process described for building the $\{\mtx{U}_{\alpha}\}_{\alpha}$ matrices.
We first recall that the difference between then HODLR and the HBS frameworks is as follows:
\begin{align*}
\mbox{HODLR framework:}&\ \mbox{The columns of $\mtb{V}_{\alpha}$ need to span the rows of }\ \mtx{A}(I_{\beta},         I_{\alpha})\\
\mbox{HBS framework:}  &\ \mbox{The columns of $\mtb{V}_{\alpha}$ need to span the rows of }\ \mtx{A}(I_{\alpha}^{\rm c},I_{\alpha}).
\end{align*}
With $\mtx{P}$ again denoting a permutation matrix, we have
$$
\mtx{A}(I_{\alpha}^{\rm c},I_{\alpha}) =
\mtx{P}
\vtwo{\mtx{A}(I_{\beta},I_{\alpha})}{\mtx{A}(I_{\tau}^{\rm c},I_{\alpha})}.
$$
It follows that the role that was played by $\mtb{Y}_{\alpha}$ in the construction of $\mtb{U}_{\alpha}$
is now played by a sampling matrix $\mtb{Z}_{\alpha}$ of size $n_{\alpha} \times r$ that satisfies
$$
\mbox{The columns of $\mtb{Z}_{\alpha}$ span the rows of $\mtx{A}(I_{\tau}^{\rm c},I_{\alpha})$.}
$$

\begin{figure}
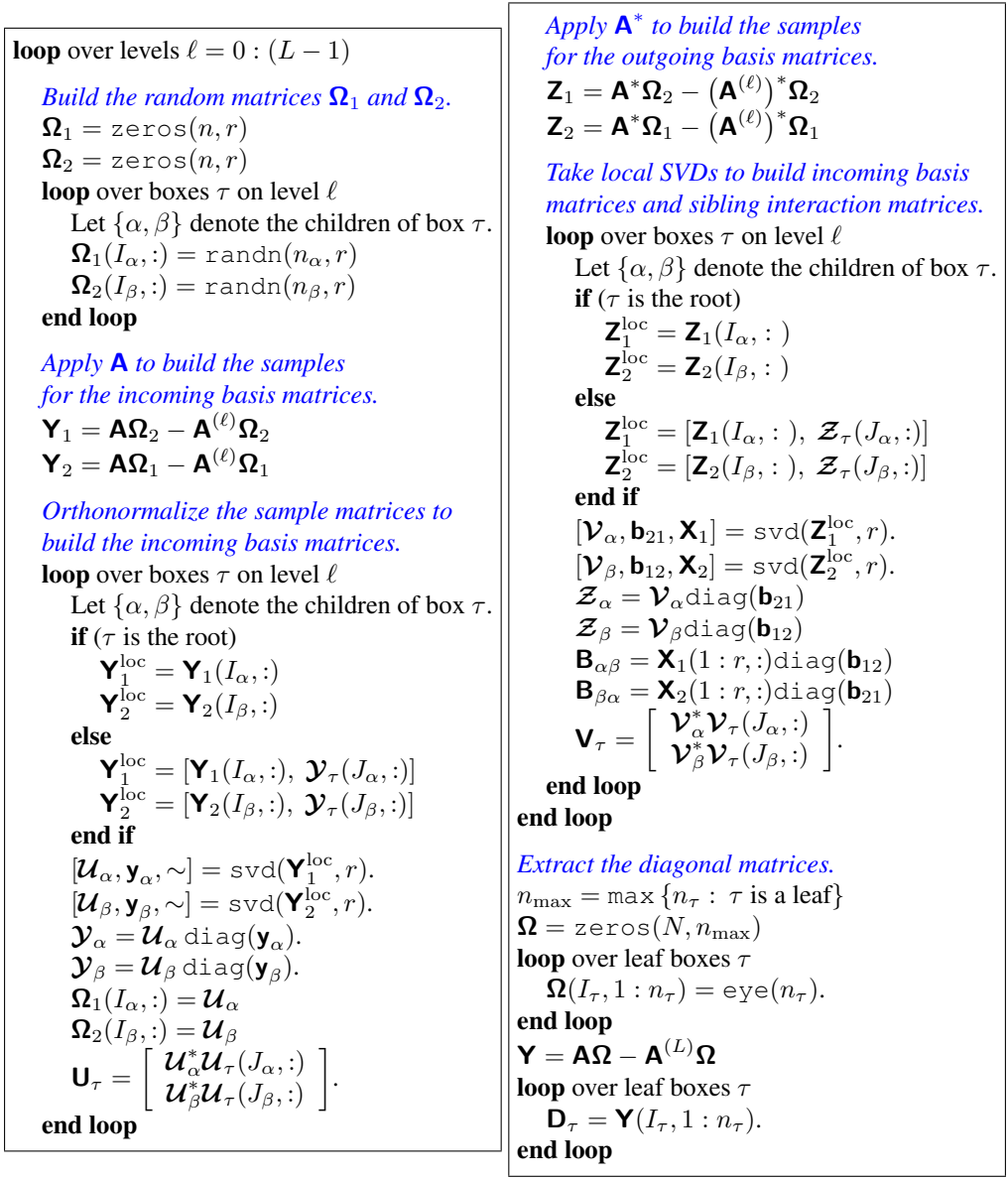

\small
\fbox{\begin{minipage}{90mm}
      \begin{tabbing}
      \hspace{3mm} \= \hspace{3mm} \= \hspace{3mm} \= \hspace{3mm}\kill
      \textbf{loop} over levels $\ell = 0:(L-1)$\\[2mm]
      \> \textit{\color{blue}Build the random matrices $\mtx{\Omega}_{1}$ and $\mtx{\Omega}_{2}$.}\\
      \> $\mtx{\Omega}_{1} = \texttt{zeros}(n,r)$\\
      \> $\mtx{\Omega}_{2} = \texttt{zeros}(n,r)$\\
      \> \textbf{loop} over boxes $\tau$ on level $\ell$\\
      \> \> Let $\{\alpha,\beta\}$ denote the children of box $\tau$.\\
      \> \> $\mtx{\Omega}_{1}(I_{\alpha},:) = \texttt{randn}(n_{\alpha},r)$\\
      \> \> $\mtx{\Omega}_{2}(I_{\beta },:) = \texttt{randn}(n_{\beta },r)$\\
      \> \textbf{end loop}\\[2mm]
      \> \textit{\color{blue}Apply $\mtx{A}$ to build the samples}\\
      \> \textit{\color{blue}for the incoming basis matrices.}\\
      \> $\mtx{Y}_{1} = \mtx{A}\mtx{\Omega}_{2} - \mtx{A}^{(\ell)}\mtx{\Omega}_{2}$\\
      \> $\mtx{Y}_{2} = \mtx{A}\mtx{\Omega}_{1} - \mtx{A}^{(\ell)}\mtx{\Omega}_{1}$\\[2mm]
      \> \textit{\color{blue}Orthonormalize the sample matrices to}\\
      \> \textit{\color{blue}build the incoming basis matrices.}\\
      \> \textbf{loop} over boxes $\tau$ on level $\ell$\\
      \> \> Let $\{\alpha,\beta\}$ denote the children of box $\tau$.\\
      \> \> \textbf{if} ($\tau$ is the root)\\
      \> \> \> $\mtx{Y}_{1}^{\rm loc} = \mtx{Y}_{1}(I_{\alpha},:)$\\
      \> \> \> $\mtx{Y}_{2}^{\rm loc} = \mtx{Y}_{2}(I_{\beta },:)$\\
      \> \> \textbf{else}\\
      \> \> \> $\mtx{Y}_{1}^{\rm loc} = [\mtx{Y}_{1}(I_{\alpha},:),\ \mtb{Y}_{\tau}(J_{\alpha},:)]$\\
      \> \> \> $\mtx{Y}_{2}^{\rm loc} = [\mtx{Y}_{2}(I_{\beta },:),\ \mtb{Y}_{\tau}(J_{\beta },:)]$\\
      \> \> \textbf{end if}\\
      \> \> $[\mtb{U}_{\alpha},\vct{y}_{\alpha},\sim] = \texttt{svd}(\mtx{Y}_{1}^{\rm loc},r)$.\\
      \> \> $[\mtb{U}_{\beta },\vct{y}_{\beta },\sim] = \texttt{svd}(\mtx{Y}_{2}^{\rm loc},r)$.\\
      \> \> $\mtb{Y}_{\alpha} = \mtb{U}_{\alpha}\,\texttt{diag}(\vct{y}_{\alpha})$.\\
      \> \> $\mtb{Y}_{\beta } = \mtb{U}_{\beta }\,\texttt{diag}(\vct{y}_{\beta })$.\\
      \> \> $\mtx{\Omega}_{1}(I_{\alpha},:) = \mtb{U}_{\alpha}$\\
      \> \> $\mtx{\Omega}_{2}(I_{\beta },:) = \mtb{U}_{\beta }$\\
      \> \> $\mtx{U}_{\tau}   = \left[\begin{array}{c} \mtb{U}_{\alpha}^{*}\mtb{U}_{\tau}(J_{\alpha},:)\\
                                                       \mtb{U}_{\beta }^{*}\mtb{U}_{\tau}(J_{\beta },:) \end{array}\right]$.\\
      \> \textbf{end loop}
      \end{tabbing}
      \end{minipage}}
\fbox{\begin{minipage}{90mm}
      \begin{tabbing}
      \hspace{3mm} \= \hspace{3mm} \= \hspace{3mm} \= \hspace{3mm}\kill
      \> \textit{\color{blue}Apply $\mtx{A}^{*}$ to build the samples}\\
      \> \textit{\color{blue}for the outgoing basis matrices.}\\
      \> $\mtx{Z}_{1} = \mtx{A}^{*}\mtx{\Omega}_{2} - \bigl(\mtx{A}^{(\ell)}\bigr)^{*}\mtx{\Omega}_{2}$\\
      \> $\mtx{Z}_{2} = \mtx{A}^{*}\mtx{\Omega}_{1} - \bigl(\mtx{A}^{(\ell)}\bigr)^{*}\mtx{\Omega}_{1}$\\[2mm]
      \> \textit{\color{blue}Take local SVDs to build incoming basis}\\
      \> \textit{\color{blue}matrices and sibling interaction matrices.}\\
      \> \textbf{loop} over boxes $\tau$ on level $\ell$\\
      \> \> Let $\{\alpha,\beta\}$ denote the children of box $\tau$.\\
      \> \> \textbf{if} ($\tau$ is the root)\\
      \> \> \> $\mtx{Z}_{1}^{\rm loc} = \mtx{Z}_{1}(I_{\alpha},\colon)$\\
      \> \> \> $\mtx{Z}_{2}^{\rm loc} = \mtx{Z}_{2}(I_{\beta },\colon)$\\
      \> \> \textbf{else}\\
      \> \> \> $\mtx{Z}_{1}^{\rm loc} = [\mtx{Z}_{1}(I_{\alpha},\colon),\ \mtb{Z}_{\tau}(J_{\alpha},:)]$\\
      \> \> \> $\mtx{Z}_{2}^{\rm loc} = [\mtx{Z}_{2}(I_{\beta },\colon),\ \mtb{Z}_{\tau}(J_{\beta },:)]$\\
      \> \> \textbf{end if}\\
      \> \> $[\mtb{V}_{\alpha},\vct{b}_{21},\mtx{X}_{1}] = \texttt{svd}(\mtx{Z}_{1}^{\rm loc},r)$.\\
      \> \> $[\mtb{V}_{\beta },\vct{b}_{12},\mtx{X}_{2}] = \texttt{svd}(\mtx{Z}_{2}^{\rm loc},r)$.\\
      \> \> $\mtb{Z}_{\alpha} = \mtb{V}_{\alpha}\texttt{diag}(\vct{b}_{21})$\\
      \> \> $\mtb{Z}_{\beta } = \mtb{V}_{\beta }\texttt{diag}(\vct{b}_{12})$\\
      \> \> $\mtx{B}_{\alpha\beta} = \mtx{X}_{1}(1:r,:)\texttt{diag}(\vct{b}_{12})$\\
      \> \> $\mtx{B}_{\beta\alpha} = \mtx{X}_{2}(1:r,:)\texttt{diag}(\vct{b}_{21})$\\
      \> \> $\mtx{V}_{\tau}   = \left[\begin{array}{c} \mtb{V}_{\alpha}^{*}\mtb{V}_{\tau}(J_{\alpha},:)\\
                                                       \mtb{V}_{\beta }^{*}\mtb{V}_{\tau}(J_{\beta },:) \end{array}\right]$.\\
      \> \textbf{end loop}\\
      \textbf{end loop}\\[2mm]
      \textit{\color{blue}Extract the diagonal matrices.}\\
      $n_{\rm max} = \texttt{max}\,\{n_{\tau}\, \colon\, \tau\mbox{ is a leaf}\}$\\
      $\mtx{\Omega} = \texttt{zeros}(N,n_{\rm max})$\\
      \textbf{loop} over leaf boxes $\tau$\\
      \> $\mtx{\Omega}(I_{\tau},1:n_{\tau}) = \texttt{eye}(n_{\tau})$.\\
      \textbf{end loop}\\
      $\mtx{Y} = \mtx{A}\mtx{\Omega} - \mtx{A}^{(L)}\mtx{\Omega}$\\
      \textbf{loop} over leaf boxes $\tau$\\
      \> $\mtx{D}_{\tau} = \mtx{Y}(I_{\tau},1:n_{\tau})$.\\
      \textbf{end loop}
      \end{tabbing}
      \end{minipage}}
\caption{A basic scheme for compressing an HBS matrix.}
\label{alg:basicHBS}
\end{figure}

The algorithm for computing the HBS representation of a matrix is now given in Figure \ref{alg:basicHBS}.

\subsection{A storage efficient scheme for compressing an HBS matrix}
\label{sec:memoryHBS}
The scheme described in Section \ref{sec:basicHBS} assumes that all ``long''
basis and spanning matrices ($\mtb{U}_{\tau}$, $\mtb{V}_{\tau}$, $\mtb{Y}_{\tau}$, $\mtb{Z}_{\tau}$)
are stored explicitly, resulting in an $O(k\,N\,\log N)$ storage requirement.
We will now demonstrate that in fact, only $O(k\,N)$ is required.

First, observe that in the HBS framework, we only need to keep at hand the long basis
matrices $\mtb{U}_{\tau}$ and $\mtb{V}_{\tau}$ for nodes $\tau$ on the level $\ell$ that is
currently being processed. In the HODLR compression algorithm in Figure \ref{alg:HODLR},
we needed the long basis matrices associated with nodes on coarser levels in order to
apply $\mtx{A}^{(\ell)}$, but in the HBS framework, all we need in order to apply $\mtx{A}^{(\ell)}$
is the long basis matrices $\{\mtb{U}_{\tau},\,\mtb{V}_{\tau}\}$ \textit{on the level
currently being processed}, and then only the short basis matrices $\mtx{U}_{\tau}$ and
$\mtx{V}_{\tau}$ for any box $\tau$ on a level coarser than $\ell$,
\textit{cf.}~the algorithm in Figure \ref{alg:apply_A_partial}.

Next, observe that the long ``spanning'' matrices $\mtb{Y}_{\tau}$ and $\mtb{Z}_{\tau}$
that were introduced in Section \ref{sec:basicHBS} do not need to be stored explicitly
either. The reason is that these matrices can be expressed in terms of the long basis
matrices $\mtb{U}_{\tau}$ and $\mtb{V}_{\tau}$. In fact, in the algorithm in Figure \ref{alg:HODLR},
we \textit{compute} $\mtb{Y}_{\tau}$ and $\mtb{Z}_{\tau}$ via the relations
$$
\mtb{Y}_{\tau} = \mtb{U}_{\tau}\,\texttt{diag}(\vct{y}_{\tau})
\qquad\mbox{and}\qquad
\mtb{Z}_{\tau} = \mtb{V}_{\tau}\,\texttt{diag}(\vct{z}_{\tau}).
$$
Since the long basis matrices $\mtb{U}_{\tau}$ and $\mtb{V}_{\tau}$ are available
during the processing of level $\ell$, we only need to store the
short vectors $\vct{y}_{\tau}$ and $\vct{z}_{\tau}$, and can then construct
$\mtb{Y}_{\tau}$ and $\mtb{Z}_{\tau}$ when they are actually needed.

\begin{figure}
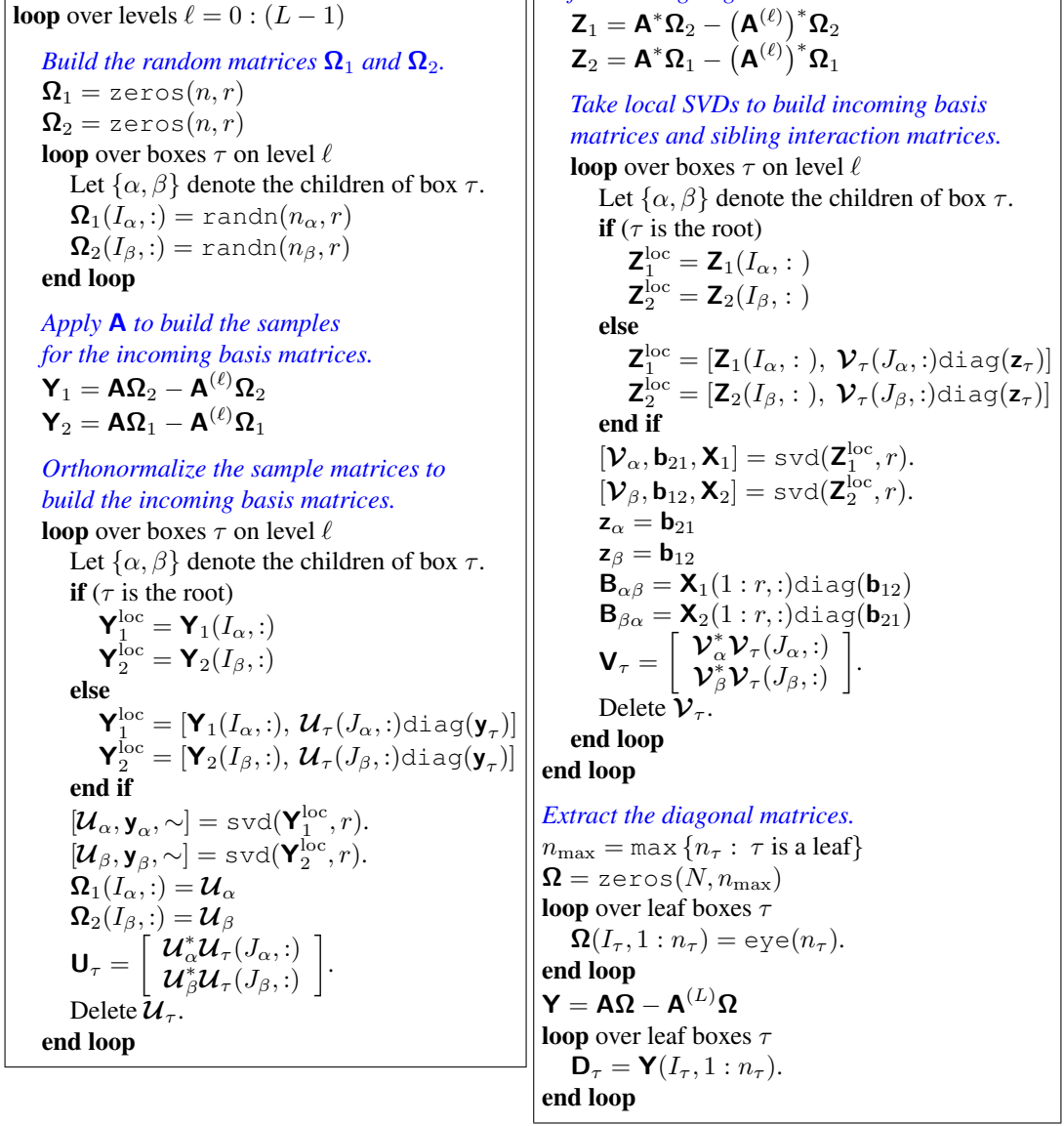

\small
\fbox{\begin{minipage}{90mm}
      \begin{tabbing}
      \hspace{3mm} \= \hspace{3mm} \= \hspace{3mm} \= \hspace{3mm}\kill
      \textbf{loop} over levels $\ell = 0:(L-1)$\\[2mm]
      \> \textit{\color{blue}Build the random matrices $\mtx{\Omega}_{1}$ and $\mtx{\Omega}_{2}$.}\\
      \> $\mtx{\Omega}_{1} = \texttt{zeros}(n,r)$\\
      \> $\mtx{\Omega}_{2} = \texttt{zeros}(n,r)$\\
      \> \textbf{loop} over boxes $\tau$ on level $\ell$\\
      \> \> Let $\{\alpha,\beta\}$ denote the children of box $\tau$.\\
      \> \> $\mtx{\Omega}_{1}(I_{\alpha},:) = \texttt{randn}(n_{\alpha},r)$\\
      \> \> $\mtx{\Omega}_{2}(I_{\beta },:) = \texttt{randn}(n_{\beta },r)$\\
      \> \textbf{end loop}\\[2mm]
      \> \textit{\color{blue}Apply $\mtx{A}$ to build the samples}\\
      \> \textit{\color{blue}for the incoming basis matrices.}\\
      \> $\mtx{Y}_{1} = \mtx{A}\mtx{\Omega}_{2} - \mtx{A}^{(\ell)}\mtx{\Omega}_{2}$\\
      \> $\mtx{Y}_{2} = \mtx{A}\mtx{\Omega}_{1} - \mtx{A}^{(\ell)}\mtx{\Omega}_{1}$\\[2mm]
      \> \textit{\color{blue}Orthonormalize the sample matrices to}\\
      \> \textit{\color{blue}build the incoming basis matrices.}\\
      \> \textbf{loop} over boxes $\tau$ on level $\ell$\\
      \> \> Let $\{\alpha,\beta\}$ denote the children of box $\tau$.\\
      \> \> \textbf{if} ($\tau$ is the root)\\
      \> \> \> $\mtx{Y}_{1}^{\rm loc} = \mtx{Y}_{1}(I_{\alpha},:)$\\
      \> \> \> $\mtx{Y}_{2}^{\rm loc} = \mtx{Y}_{2}(I_{\beta },:)$\\
      \> \> \textbf{else}\\
      \> \> \> $\mtx{Y}_{1}^{\rm loc} = [\mtx{Y}_{1}(I_{\alpha},:),\ \mtb{U}_{\tau}(J_{\alpha},:)\texttt{diag}(\vct{y}_{\tau})]$\\
      \> \> \> $\mtx{Y}_{2}^{\rm loc} = [\mtx{Y}_{2}(I_{\beta },:),\ \mtb{U}_{\tau}(J_{\beta },:)\texttt{diag}(\vct{y}_{\tau})]$\\
      \> \> \textbf{end if}\\
      \> \> $[\mtb{U}_{\alpha},\vct{y}_{\alpha},\sim] = \texttt{svd}(\mtx{Y}_{1}^{\rm loc},r)$.\\
      \> \> $[\mtb{U}_{\beta },\vct{y}_{\beta },\sim] = \texttt{svd}(\mtx{Y}_{2}^{\rm loc},r)$.\\
      \> \> $\mtx{\Omega}_{1}(I_{\alpha},:) = \mtb{U}_{\alpha}$\\
      \> \> $\mtx{\Omega}_{2}(I_{\beta },:) = \mtb{U}_{\beta }$\\
      \> \> $\mtx{U}_{\tau}   = \left[\begin{array}{c} \mtb{U}_{\alpha}^{*}\mtb{U}_{\tau}(J_{\alpha},:)\\
                                                       \mtb{U}_{\beta }^{*}\mtb{U}_{\tau}(J_{\beta },:) \end{array}\right]$.\\
      \> \> Delete $\mtb{U}_{\tau}$.\\
      \> \textbf{end loop}
      \end{tabbing}
      \end{minipage}}
\fbox{\begin{minipage}{90mm}
      \begin{tabbing}
      \hspace{3mm} \= \hspace{3mm} \= \hspace{3mm} \= \hspace{3mm}\kill
      \> \textit{\color{blue}Apply $\mtx{A}^{*}$ to build the samples}\\
      \> \textit{\color{blue}for the outgoing basis matrices.}\\
      \> $\mtx{Z}_{1} = \mtx{A}^{*}\mtx{\Omega}_{2} - \bigl(\mtx{A}^{(\ell)}\bigr)^{*}\mtx{\Omega}_{2}$\\
      \> $\mtx{Z}_{2} = \mtx{A}^{*}\mtx{\Omega}_{1} - \bigl(\mtx{A}^{(\ell)}\bigr)^{*}\mtx{\Omega}_{1}$\\[2mm]
      \> \textit{\color{blue}Take local SVDs to build incoming basis}\\
      \> \textit{\color{blue}matrices and sibling interaction matrices.}\\
      \> \textbf{loop} over boxes $\tau$ on level $\ell$\\
      \> \> Let $\{\alpha,\beta\}$ denote the children of box $\tau$.\\
      \> \> \textbf{if} ($\tau$ is the root)\\
      \> \> \> $\mtx{Z}_{1}^{\rm loc} = \mtx{Z}_{1}(I_{\alpha},\colon)$\\
      \> \> \> $\mtx{Z}_{2}^{\rm loc} = \mtx{Z}_{2}(I_{\beta },\colon)$\\
      \> \> \textbf{else}\\
      \> \> \> $\mtx{Z}_{1}^{\rm loc} = [\mtx{Z}_{1}(I_{\alpha},\colon),\ \mtb{V}_{\tau}(J_{\alpha},:)\texttt{diag}(\vct{z}_{\tau})]$\\
      \> \> \> $\mtx{Z}_{2}^{\rm loc} = [\mtx{Z}_{2}(I_{\beta },\colon),\ \mtb{V}_{\tau}(J_{\beta },:)\texttt{diag}(\vct{z}_{\tau})]$\\
      \> \> \textbf{end if}\\
      \> \> $[\mtb{V}_{\alpha},\vct{b}_{21},\mtx{X}_{1}] = \texttt{svd}(\mtx{Z}_{1}^{\rm loc},r)$.\\
      \> \> $[\mtb{V}_{\beta },\vct{b}_{12},\mtx{X}_{2}] = \texttt{svd}(\mtx{Z}_{2}^{\rm loc},r)$.\\
      \> \> $\vct{z}_{\alpha} = \vct{b}_{21}$\\
      \> \> $\vct{z}_{\beta } = \vct{b}_{12}$\\
      \> \> $\mtx{B}_{\alpha\beta} = \mtx{X}_{1}(1:r,:)\texttt{diag}(\vct{b}_{12})$\\
      \> \> $\mtx{B}_{\beta\alpha} = \mtx{X}_{2}(1:r,:)\texttt{diag}(\vct{b}_{21})$\\
      \> \> $\mtx{V}_{\tau}   = \left[\begin{array}{c} \mtb{V}_{\alpha}^{*}\mtb{V}_{\tau}(J_{\alpha},:)\\
                                                       \mtb{V}_{\beta }^{*}\mtb{V}_{\tau}(J_{\beta },:) \end{array}\right]$.\\
      \> \> Delete $\mtb{V}_{\tau}$.\\
      \> \textbf{end loop}\\
      \textbf{end loop}\\[2mm]
      \textit{\color{blue}Extract the diagonal matrices.}\\
      $n_{\rm max} = \texttt{max}\,\{n_{\tau}\, \colon\, \tau\mbox{ is a leaf}\}$\\
      $\mtx{\Omega} = \texttt{zeros}(N,n_{\rm max})$\\
      \textbf{loop} over leaf boxes $\tau$\\
      \> $\mtx{\Omega}(I_{\tau},1:n_{\tau}) = \texttt{eye}(n_{\tau})$.\\
      \textbf{end loop}\\
      $\mtx{Y} = \mtx{A}\mtx{\Omega} - \mtx{A}^{(L)}\mtx{\Omega}$\\
      \textbf{loop} over leaf boxes $\tau$\\
      \> $\mtx{D}_{\tau} = \mtx{Y}(I_{\tau},1:n_{\tau})$.\\
      \textbf{end loop}
      \end{tabbing}
      \end{minipage}}
\caption{A storage efficient algorithm for compressing an HBS matrix.}
\label{alg:memoryHBS}
\end{figure}

The memory efficient algorithm resulting from exploiting the observations described
in this section is summarized in Figure \ref{alg:memoryHBS}.

\subsection{Adaptive rank determination and conversion to the HBS-ID format}
\label{sec:HBSID}
The schemes presented in Sections \ref{sec:basicHBS} and \ref{sec:memoryHBS} do not adaptively
determine the ranks of the off-diagonal blocks being compressed. Instead, every block is
factored using a preset uniform rank $\ell$ that must be picked to be larger than any
actual numerical rank encountered. It is possible to incorporate adaptive rank determination
in to the scheme, but we found it easier to perform this step in a second sweep that
travels through the tree in the opposite direction: from smaller boxes to larger. In this
second sweep, we also convert the standard HBS format to the HBS-ID format, which leads to
a slight improvement in storage requirements, and improves interpretability of the
sibling interaction matrices, as discussed in Section \ref{sec:introHBS}
and Definition \ref{def:HBSID}.

The conversion to the HBS-ID is a ``post-processing'' step, so in what follows, we assume
that the compression algorithm in Figure \ref{alg:memoryHBS} has already been executed
so that both the HBS basis matrices $\mtb{U}_{\tau}$ and $\mtb{V}_{\tau}$, and the
sample matrix $\mtb{Y}_{\tau}$ and $\mtb{Z}_{\tau}$ are available for every node
(these are stored \textit{implicitly} in terms of the short basis matrices
$\mtx{U}_{\tau}$ and $\mtx{V}_{\tau}$, as described in Section \ref{sec:memoryHBS}).

The first step is to sweep over all leaves $\tau$ in the tree. For each leaf, we now have
available spanning matrices $\mtb{Y}_{\tau}$ and $\mtb{Z}_{\tau}$ whose columns span
the columns of $\mtx{A}(I_{\tau},I_{\tau}^{\rm c})$ and $\mtx{A}(I_{\tau}^{\rm c},I_{\tau})^{*}$,
respectively. In order to find a set of spanning rows $\tilde{I}_{\tau}^{\rm in}$ of $\mtx{A}(I_{\tau},I_{\tau}^{\rm c})$
and a set of spanning columns $\tilde{I}_{\tau}^{\rm out}$ of $\mtx{A}(I_{\tau}^{\rm c},I_{\tau})^{*}$,
all we need to do is to compute interpolatory decompositions (IDs) of the small matrices
$\mtb{Y}_{\tau}$ and $\mtb{Z}_{\tau}$, cf.~Remark \ref{remark:spanningskel},
\begin{equation}
\label{eq:buildTJ}
[\mtx{T}_{\rm in },J_{\rm in }] = \texttt{id}(\mtb{Y}_{\tau}^{*},\varepsilon),
\qquad\mbox{and}\qquad
[\mtx{T}_{\rm out},J_{\rm out}] = \texttt{id}(\mtb{Z}_{\tau}^{*},\varepsilon).
\end{equation}
In equation (\ref{eq:buildTJ}), we give the computational tolerance $\varepsilon$ as an
input parameter. This reveals the ``true'' $\varepsilon$-ranks $k_{\rm in}$ and $k_{\rm out}$.
Then the skeleton index vectors for $\tau$ are given by
$$
\tilde{I}_{\tau}^{\rm in } = I_{\tau}(J_{\rm in }(1:k_{\rm in })),
\qquad\mbox{and}\qquad
\tilde{I}_{\tau}^{\rm out} = I_{\tau}(J_{\rm out}(1:k_{\rm out})).
$$
Now define the \textit{subsampled} basis matrices $\mtx{U}_{\tau}^{\rm samp}$ and $\mtx{V}_{\tau}^{\rm samp}$ via
\begin{equation}
\label{eq:ben0}
\mtx{U}_{\tau}^{\rm samp} = \mtx{U}(J_{\rm in }(1:k_{\rm in }),:),
\qquad\mbox{and}\qquad
\mtx{V}_{\tau}^{\rm samp} = \mtx{V}(J_{\rm out}(1:k_{\rm out}),:).
\end{equation}

Once all leaves have been processed in this manner, we can determine the sibling interaction matrices
in the HBS-ID representation. Let $\{\alpha,\beta\}$ denote a sibling pair consisting
of two leaves. First observe that, by definition,
\begin{equation}
\label{eq:ben1}
\mtx{B}^{\rm skel}_{\alpha,\beta} = \mtx{A}(\tilde{I}_{\alpha}^{\rm in},\tilde{I}_{\beta}^{\rm out}).
\end{equation}
Next, recall that
\begin{equation}
\label{eq:ben2}
\mtx{A}(I_{\alpha},I_{\beta}) = \mtx{U}_{\alpha}\mtx{B}_{\alpha,\beta}\mtx{V}_{\beta}^{*}.
\end{equation}
Combining (\ref{eq:ben0}), (\ref{eq:ben1}), and (\ref{eq:ben2}), we find that
$$
\mtx{B}^{\rm skel}_{\alpha,\beta} = \mtx{U}_{\alpha}^{\rm samp}\,\mtx{B}_{\alpha,\beta}\,(\mtx{V}_{\beta}^{\rm samp})^{*}.
$$

Once all leaves have been processed, we next proceed to the parent nodes. We do this
by transversing the tree in the other direction, going from smaller to larger boxes.
When a box $\tau$ is processed, its children $\{\alpha,\beta\}$ have already been
processed. The key observation is now that that if we set
$\hat{I}_{\tau}^{\rm in } = \tilde{I}_{\alpha}^{\rm in } \cup \tilde{I}_{\beta}^{\rm in }$ and
$\hat{I}_{\tau}^{\rm out} = \tilde{I}_{\alpha}^{\rm out} \cup \tilde{I}_{\beta}^{\rm out}$,
then these index vectors form skeletons for $\tau$. These skeletons are inefficient, but
by simply compressing the corresponding rows of $\mtb{Y}_{\tau}$ and $\mtb{Z}_{\tau}$,
we can build the skeletons and the interpolation matrices associated with $\tau$. Due
to the self-similarity between levels in the HBS representation, the compression is
entirely analogous to the compression of a leaf, with the only difference that the
role played by the basis matrices $\mtb{U}_{\tau}$ and $\mtb{V}_{\tau}$ for a leaf,
are now played by the sub-sampled matrices $\mtx{U}_{\rm tmp}$ and $\mtx{V}_{\rm tmp}$
which represent the restriction of $\mtb{U}_{\tau}$ and $\mtb{V}_{\tau}$ to the index
rows and columns indicated by the index vectors $\hat{I}_{\tau}^{\rm in}$ and
$\hat{I}_{\tau}^{\rm out}$, respectively. The entire process is summarized in Figure \ref{alg:HBSID}.

\begin{figure}
\small
\fbox{\begin{minipage}{90mm}
      \begin{tabbing}
      \hspace{3mm} \= \hspace{3mm} \= \hspace{3mm} \= \hspace{3mm}\kill
      Execute the compression algorithm described in Figure \ref{alg:memoryHBS}.\\
      \textbf{loop} over leaves $\tau$ \\
      \> $\mtx{Y}_{\rm tmp} = \mtx{U}_{\tau}\texttt{diag}(\vct{y}_{\tau})$\\
      \> $\mtx{Z}_{\rm tmp} = \mtx{V}_{\tau}\texttt{diag}(\vct{z}_{\tau})$\\
      \> $[\mtx{T}_{\rm in },J_{\rm in }] = \texttt{id}(\mtx{Y}_{\rm tmp}^{*},\varepsilon)$\\
      \> $[\mtx{T}_{\rm out},J_{\rm out}] = \texttt{id}(\mtx{Z}_{\rm tmp}^{*},\varepsilon)$\\
      \> $\mtx{U}_{\tau}^{\rm samp} = \mtx{U}_{\tau}(J_{\rm in }(1:k_{\rm in}),:)$.\\
      \> $\mtx{V}_{\tau}^{\rm samp} = \mtx{V}_{\tau}(J_{\rm out}(1:k_{\rm out}),:)$.\\
      \textbf{end loop}
      \end{tabbing}
      \end{minipage}}
\fbox{\begin{minipage}{90mm}
      \begin{tabbing}
      \hspace{3mm} \= \hspace{3mm} \= \hspace{3mm} \= \hspace{3mm}\kill
      \textbf{loop} over levels $\ell = (L-1):(-1):1$\\
      \> \textbf{loop} over boxes $\tau$ on level $\ell$\\
      \> \> Let $\{\alpha,\beta\}$ denote the children of box $\tau$.\\
      \> \> $\mtx{U}_{\rm tmp} = \left[\begin{array}{cc}
                                       \mtx{U}_{\alpha}^{\rm samp} & \mtx{0} \\
                                       \mtx{0} & \mtx{U}_{\beta}^{\rm samp}
                                       \end{array}\right]
                                       \mtx{U}_{\tau}$.\\
      \> \> $\mtx{V}_{\rm tmp} = \left[\begin{array}{cc}
                                       \mtx{V}_{\alpha}^{\rm samp} & \mtx{0} \\
                                       \mtx{0} & \mtx{V}_{\beta}^{\rm samp}
                                       \end{array}\right]
                                       \mtx{V}_{\tau}$.\\
      \> \> $\mtx{Y}_{\rm tmp} = \mtx{U}_{\rm tmp}\texttt{diag}(\vct{y}_{\tau})$\\
      \> \> $\mtx{Z}_{\rm tmp} = \mtx{V}_{\rm tmp}\texttt{diag}(\vct{z}_{\tau})$\\
      \> \> $[\mtx{T}_{\rm in },J_{\rm in }] = \texttt{id}(\mtx{Y}_{\rm tmp}^{*},\varepsilon)$\\
      \> \> $[\mtx{T}_{\rm out},J_{\rm out}] = \texttt{id}(\mtx{Z}_{\rm tmp}^{*},\varepsilon)$\\
      \> \> $\mtx{U}_{\tau}^{\rm samp} = \mtx{U}_{\rm tmp}(J_{\rm in }(1:k_{\rm in }),:)$.\\
      \> \> $\mtx{V}_{\tau}^{\rm samp} = \mtx{V}_{\rm tmp}(J_{\rm out}(1:k_{\rm out}),:)$.\\
      \> \> $\mtx{B}_{\alpha,\beta}^{\rm skel} = \mtx{U}_{\alpha}^{\rm samp}\mtx{B}_{\alpha,\beta}(\mtx{V}_{\beta}^{\rm samp})^{*}$.\\
      \> \> $\mtx{B}_{\beta,\alpha}^{\rm skel} = \mtx{U}_{\beta}^{\rm samp}\mtx{B}_{\beta,\alpha}(\mtx{V}_{\alpha}^{\rm samp})^{*}$.\\
      \> \textbf{end loop}\\
      \textbf{end loop}\\
      Let $\{\alpha,\beta\}$ denote the children of the root.\\
      $\mtx{B}_{\alpha,\beta}^{\rm skel} = \mtx{U}_{\alpha}^{\rm samp}\mtx{B}_{\alpha,\beta}(\mtx{V}_{\beta}^{\rm samp})^{*}$.\\
      $\mtx{B}_{\beta,\alpha}^{\rm skel} = \mtx{U}_{\beta}^{\rm samp}\mtx{B}_{\beta,\alpha}(\mtx{V}_{\alpha}^{\rm samp})^{*}$.
      \end{tabbing}
      \end{minipage}}
\caption{An algorithm for computing the HBS-ID representation of a given matrix $\mtx{A}$. This algorithm
adaptively determines the ranks of the off-diagonal blocks of $\mtx{A}$.}
\label{alg:HBSID}
\end{figure}

\subsection{Asymptotic complexity}
\label{sec:HBSasympt}
The asymptotic complexity for the HBSID algorithm is very similar to that for the HODLR algorithm.
The key difference is that since $\mtx{A}^{(\ell)}$ is now applied using nested basis, we find
$$
T_{\mtx{A}^{(\ell)}}
\sim T_{\rm flop} \times \left(2^{\ell}\,k\,\frac{N}{2^{\ell}} + \sum_{j=0}^{\ell-1} 2^{j}\,k^{2} \right)
\sim T_{\rm flop} \times k\,N.
$$
In other words, the asymptotic complexity of applying $T_{\mtx{A}^{(\ell)}}$ is now less by a factor of $O(\ell)$.
In consequence,
\begin{equation}
\label{eq:T_HBSID}
T_{\rm compress} \sim T_{\rm mult} \times k\,\log N + T_{\rm flop} \times k^{2}\,N\,\log N.
\end{equation}

\section{Numerical experiments}
\label{sec:num}

In this section, we present results from numerical experiments that substantiate
claims on asymptotic complexity made in sections \ref{sec:HODLRasympt} and \ref{sec:HBSasympt},
and demonstrate that the practical execution time is very competitive (in other words, that
the scaling constants suppressed in the asymptotic analysis are moderate). We investigate
four different test problems:
In Section \ref{sec:BIE} we apply the randomized compression schemes to a discretized boundary
integral operator for which other compression techniques are already available. This allows us
to benchmark the new algorithms and verify their accuracy.
In Section \ref{sec:DfN} we demonstrate how the proposed scheme can be used to form a compressed
representation of a product of two compressed matrices, thus demonstrated a way of circumventing
the need for a complex and time-consuming structured matrix-matrix multiplication.
In Section \ref{sec:FMM} we apply the schemes to a potential evaluation problem where we use
the Fast Multipole Method (FMM) as the ``black-box'' matrix-vector multiplication scheme (note that the
data sparse format implicit in the FMM is much more cumbersome to invert than the HODLR and HBSID
formats).
In Section \ref{sec:nesteddiss} we apply the scheme to compress large dense matrices that
arise in the classical ``nested dissection'' or ``multifrontal'' direct solvers for
the sparse matrices arising from finite element or finite difference discretization of
elliptic PDEs. The ability to efficiently manipulate such matrices allows for the
construction of $O(N)$ complexity direct solvers for the associated sparse linear systems.

For each of the four test problems, we compare two different techniques for computing a
data-sparse representation of $\mtx{A}$: (1) The randomized technique for computing an
HODLR-representation in Figure \ref{alg:HODLR}. (2) The randomized technique for
computing an HBSID-representation in Figure \ref{alg:HBSID}. For each technique, we
report the following quantities:
\begin{center}
\begin{tabular}{lll}
$N$ &\mbox{}\hspace{1mm}\mbox{}& The number of DOFs (so that $\mtx{A}$ is of size $N\times N$).\\
$\ell$ && The number of random vectors used at each level ($\ell$ must be larger than the maximal rank).\\
$k$ && The largest rank encountered in the compression.\\
$N_{\rm matvec}$ && The number of applications of $\mtx{A}$ required (so that $N_{\rm matvec} = (L+1)\times \ell \sim \log(N)\times \ell$).\\
$T_{\rm compress}$ && The time required for compression (in seconds).\\
$T_{\rm net}$ && The time required for compression, excluding time to apply $\mtx{A}$ and $\mtx{A}^{*}$ (in seconds).\\
$T_{\rm app}$ && The time required for applying the compressed matrix to a vector (in seconds).\\
$M$ && The amount of memory required to store $\mtx{A}$ (in MB).
\end{tabular}
\end{center}
The reason that we report the time $T_{\rm net}$ (that does not count time spent in the black-box
matrix-vector multiplier) is to validate our claims (\ref{eq:T_HODLR}) and (\ref{eq:T_HBSID})
regarding the asymptotic complexity of the method. To summarize, our predictions are, for the
HODLR algorithm
$$
T_{\rm net} \sim N\,(\log N)^{2},
\qquad
T_{\rm app} \sim N\,\log N,
\qquad
M \sim N\,\log(N),
$$
and for the HBSID algorithm
$$
T_{\rm net} \sim N\,\log N,
\qquad
T_{\rm app} \sim N,
\qquad
M \sim N.
$$
In addition to the timings, we computed a randomized estimate $E$ of the compression error, computed
as follows: We drew ten vectors $\{\vct{\omega}_{i}\}_{i=1}^{10}$ of unit length from a uniform
distribution on the unit sphere in $\mathbb{R}^{N}$. Then $E$ is defined via
$$
E = \max_{1 \leq i \leq 10}
\frac{||\mtx{A}\vct{\omega}_{i} - \mtx{A}_{\rm compressed}\vct{\omega}_{i}||}
     {||\mtx{A}\vct{\omega}_{i}||}.
$$

\subsection{Compressing a Boundary Integral Equation}
\label{sec:BIE}
\begin{figure}
\includegraphics[height=30mm]{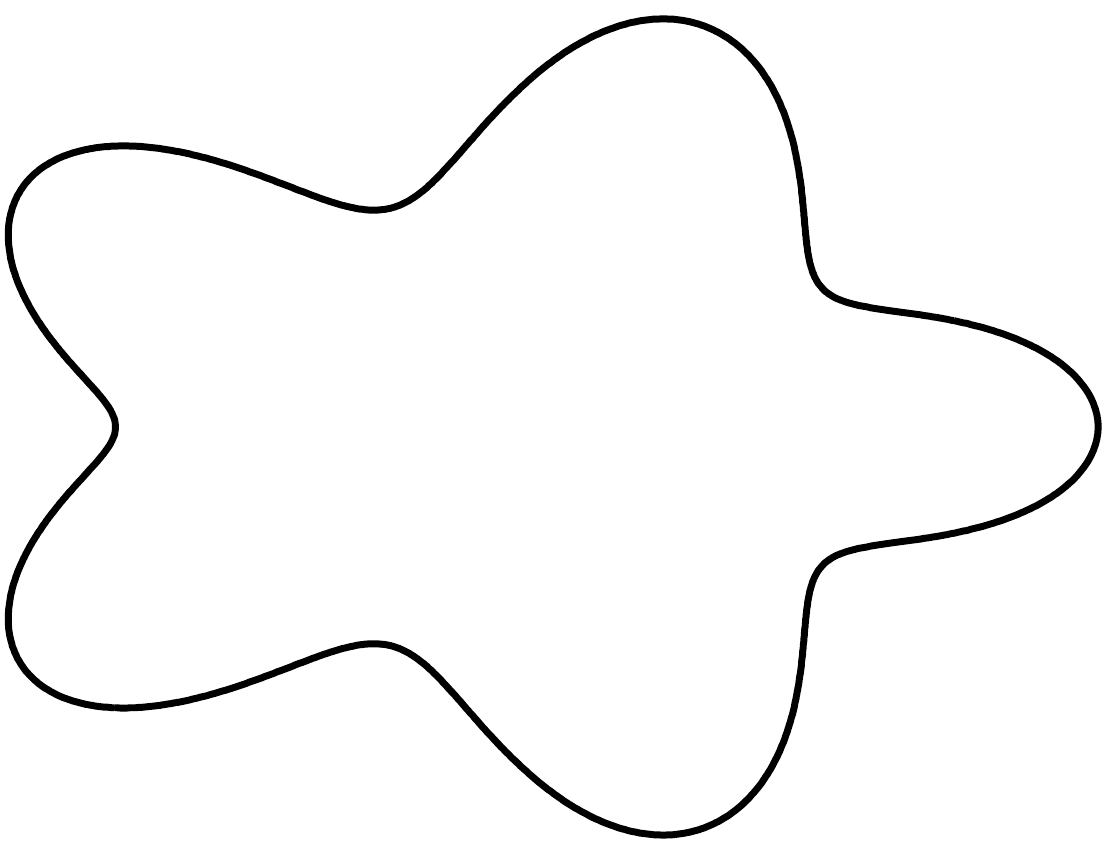}
\caption{Contour $\Gamma$ on which the BIE (\ref{eq:BIE}) in Section \ref{sec:BIE} is defined.}
\label{fig:BIE}
\end{figure}

Our first numerical example concerns compression of a discretized
version of the Boundary Integral Equation (BIE)
\begin{equation}
\label{eq:BIE}
\frac{1}{2}q(\pvct{x}) +
\int_{\Gamma}
\frac{(\pvct{x} - \pvct{y})\cdot\pvct{n}(\pvct{y})}
     {4\pi|\pvct{x} - \pvct{y}|^{2}}\,q(\pvct{y})\,ds(\pvct{y}) = f(\pvct{x}),
\qquad\pvct{x} \in \Gamma,
\end{equation}
where $\Gamma$ is the simple curve shown in Figure \ref{fig:BIE}, and where
$\pvct{n}(\pvct{y})$ is the outwards pointing unit normal of $\Gamma$ at
$\pvct{y}$. The BIE (\ref{eq:BIE}) is a standard integral equation formulation
of the Laplace equation with boundary condition $f$ on the domain interior
to $\Gamma$.\footnote{Verify!} We discretize the BIE (\ref{eq:BIE}) using the
Nystr\"om method on $N$ equispaced points on $\Gamma$, with the Trapezoidal
rule as the quadrature. Note that the kernel in (\ref{eq:BIE}) is smooth, so
the Trapezoidal rule has exponential convergence. This problem is slightly
artificial in that only about 200 points \footnote{Check.} are needed to attain
full double precision accuracy in the discretization. We include it for bench-marking
purposes to verify the scaling of the proposed method.

To be precise, the matrix $\mtx{A}$ used in this numerical experiment is
itself an HBS representation of the matrix resulting from discretization
of (\ref{eq:BIE}), computed using the technique of \cite{2005_martinsson_fastdirect}.
To minimize the risk of spurious effects due to both the ``exact'' and the
computed $\mtx{A}$ being HBS representations, we used a much higher
precision in computing the ``exact'' $\mtx{A}$, and also a shifted tree
structure to avoid having the compressed blocks of our reference $\mtx{A}$
align with the compressed blocks constructed by the randomzed sampling algorithms.

\begin{figure}
\includegraphics[width=\textwidth]{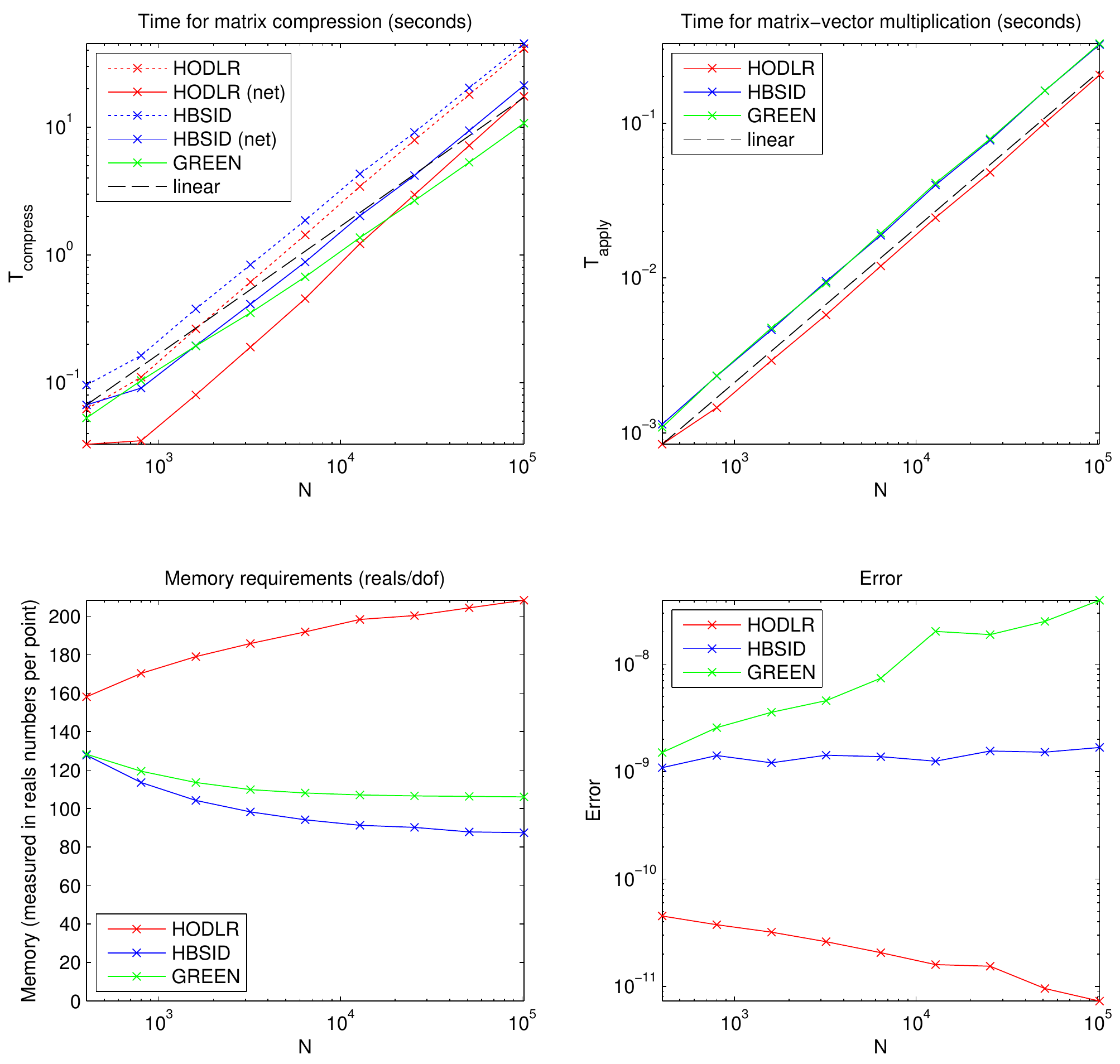}
\caption{Visualization of the results presented in Tables \ref{table:BIE_S}, \ref{table:BIE_HBS}, and \ref{table:BIE_Green},
pertaining to the example in Section \ref{sec:BIE}.}
\label{fig:BIE_results}
\end{figure}

For this experiment, we bench-mark the new compression algorithms by comparing
their speed, accuracy, and memory requirements to the compression technique
based on potential theory described in \cite{2005_martinsson_fastdirect},
run at the same precision as the randomized compression scheme. The results are
presented in tables \ref{table:BIE_S}, \ref{table:BIE_HBS}, and \ref{table:BIE_Green},
and summarized in Figure \ref{fig:BIE_results}.

\begin{table}
\begin{center}
\begin{tabular}{|r|r|r|r|r|r|r|r|r|r|r|}\hline
$N$ & $N_{\rm matvec}$ & $T_{\rm compress}$ & $T_{\rm net}$ & $T_{\rm app}$ & $M$ & $M/N$ & $E$ & $k$\\
&& (sec)& (sec) & (sec) & (MB) & (reals) & & \\ \hline
   400 &   3 x 35 &     0.065 &      0.035 &     0.001 &    0.5 & 157.7 &  4.68e-11 &  28 \\
   800 &   4 x 35 &     0.107 &      0.033 &     0.001 &    1.0 & 169.9 &  3.67e-11 &  28 \\
  1600 &   5 x 35 &     0.267 &      0.081 &     0.003 &    2.2 & 179.1 &  3.09e-11 &  28 \\
  3200 &   6 x 35 &     0.617 &      0.188 &     0.006 &    4.5 & 186.2 &  2.57e-11 &  28 \\
  6400 &   7 x 35 &     1.458 &      0.456 &     0.011 &    9.5 & 194.0 &  1.99e-11 &  29 \\
 12800 &   8 x 35 &     3.453 &      1.198 &     0.024 &   19.6 & 200.5 &  1.80e-11 &  30 \\
 25600 &   9 x 35 &     8.065 &      2.998 &     0.047 &   39.7 & 203.1 &  1.74e-11 &  29 \\
 51200 &  10 x 35 &    18.362 &      7.140 &     0.099 &   82.1 & 210.3 &  1.14e-11 &  30 \\
102400 &  11 x 35 &    41.985 &     17.366 &     0.204 &  166.6 & 213.2 &  1.24e-11 &  30 \\ \hline
\end{tabular}
\end{center}
\caption{Compression to the HODLR format using Algorithm \ref{alg:HODLR} of the double layer integral equation described in
Section \ref{sec:BIE}. Here $\varepsilon = 10^{-9}$ and $\ell = 35$.}
\label{table:BIE_S}
\end{table}

\begin{table}
\begin{center}
\begin{tabular}{|r|r|r|r|r|r|r|r|r|r|r|}\hline
$N$ & $N_{\rm matvec}$ & $T_{\rm compress}$ & $T_{\rm net}$ & $T_{\rm app}$ & $M$ & $M/N$ & $E$ & $k$\\
&& (sec)& (sec) & (sec) & (MB) & (reals) & & \\ \hline
   400 &   3 x 35 &     0.098 &      0.068 &     0.001 &    0.4 & 127.7 &  1.40e-09 &  25 \\
   800 &   4 x 35 &     0.163 &      0.090 &     0.002 &    0.7 & 113.6 &  1.30e-09 &  25 \\
  1600 &   5 x 35 &     0.381 &      0.193 &     0.005 &    1.3 & 104.3 &  1.42e-09 &  24 \\
  3200 &   6 x 35 &     0.840 &      0.411 &     0.009 &    2.4 &  98.3 &  1.32e-09 &  24 \\
  6400 &   7 x 35 &     1.881 &      0.878 &     0.019 &    4.6 &  94.2 &  1.56e-09 &  24 \\
 12800 &   8 x 35 &     4.258 &      1.969 &     0.037 &    8.9 &  91.3 &  2.04e-09 &  23 \\
 25600 &   9 x 35 &     9.262 &      4.205 &     0.076 &   17.6 &  90.3 &  1.62e-09 &  23 \\
 51200 &  10 x 35 &    20.431 &      9.238 &     0.153 &   34.4 &  87.9 &  1.44e-09 &  23 \\
102400 &  11 x 35 &    45.732 &     21.039 &     0.310 &   68.3 &  87.5 &  1.61e-09 &  23 \\ \hline
\end{tabular}
\end{center}
\caption{Compression to the HBSID format using Algorithm \ref{alg:HBSID} of the double layer integral equation described in
Section \ref{sec:BIE}. Here $\varepsilon = 10^{-9}$ and $\ell = 35$.}
\label{table:BIE_HBS}
\end{table}

\begin{table}
\begin{center}
\begin{tabular}{|l|r|r|r|r|r|r|r|r|r|r|}\hline
$N$ & $T_{\rm compress}$ & $t_{\rm app}$ & $M$ & $M/N$ & $E$ & $k$\\
& (sec)& (sec) & (MB) & (reals) & & \\ \hline
   400 &     0.053 &     0.001 &    0.4 & 128.3 &  1.51e-09 &  38 \\
   800 &     0.104 &     0.002 &    0.7 & 119.5 &  2.57e-09 &  37 \\
  1600 &     0.194 &     0.005 &    1.4 & 113.6 &  3.57e-09 &  35 \\
  3200 &     0.352 &     0.009 &    2.7 & 109.9 &  4.58e-09 &  34 \\
  6400 &     0.674 &     0.019 &    5.3 & 108.2 &  7.42e-09 &  33 \\
 12800 &     1.368 &     0.041 &   10.5 & 107.1 &  2.02e-08 &  31 \\
 25600 &     2.664 &     0.079 &   20.8 & 106.7 &  1.89e-08 &  29 \\
 51200 &     5.308 &     0.162 &   41.5 & 106.3 &  2.51e-08 &  28 \\
102400 &    10.758 &     0.327 &   82.9 & 106.1 &  3.95e-08 &  25 \\ \hline
\end{tabular}
\end{center}
\caption{Compression to the HBSID format of the double layer integral equation described in
Section \ref{sec:BIE}, using the technique based on potential theory of \cite{2005_martinsson_fastdirect}
with $\varepsilon = 10^{-9}$.}
\label{table:BIE_Green}
\end{table}

\subsection{Operator multiplication}
\label{sec:DfN}
We next apply the proposed scheme to compute the Neumann-to-Dirichlet operator for the
boundary value problem
\begin{equation}
\label{eq:BVP}
\left\{\begin{aligned}
-\Delta u(\pvct{x}) =&\ 0\qquad&\pvct{x} \in \Omega,\\
\partial_{\pvct{n}}u(\pvct{x}) =&\ g(\pvct{x})\qquad&\pvct{x} \in \Gamma,
\end{aligned}\right.
\end{equation}
where $\Gamma$ is again the contour shown in Figure \ref{fig:BIE}, where $\Omega$
is the domain exterior to $\Gamma$, and where $\pvct{n}$ is the unit normal vector
pointing in the outwards direction from $\Gamma$. With $u$ the solution of (\ref{eq:BVP}), let
$f$ denote the restriction of $u$ to $\Gamma$, and let $T$ denote the linear operator
$
T\,\colon\,g \mapsto f,
$
known as the \textit{Neumann-to-Dirichlet (NtD) operator}. It is well-known (see Remark \ref{remark:DfN})
that $T$ can be built explicitly as the product
\begin{equation}
\label{eq:DfN}
T = S\,\left(\frac{1}{2}I + D^{*}\right)^{-1},
\end{equation}
where $S$ is the \textit{single-layer operator}
$
[Sq](\pvct{x}) = \int_{\Gamma}-\frac{1}{2\pi}\,\log|\pvct{x} - \pvct{y}|\,q(\pvct{y})\,ds(\pvct{y}),
$
and where $D^{*}$ is the adjoint of the double-layer operator
$
[D^{*}q](\pvct{x}) = \int_{\Gamma}\frac{\pvct{n}(\pvct{x}) \cdot (\pvct{x} - \pvct{y})}{2\pi |\pvct{x} - \pvct{y}|^{2}}
\,q(\pvct{y})\,ds(\pvct{y}).
$
We form discrete approximations $\mtx{S}$ and $\mtx{D}^{*}$ to $S$ and $D^{*}$, and compute
$(\tfrac{1}{2}\mtx{I} + \mtx{D}^{*})^{-1}$ using the techniques in \cite{2012_martinsson_FDS_survey},
with $6^{\rm th}$ order Kapur-Rokhlin quadrature used to discretize the singular integral
operator $S$. Then we can evaluate a discrete approximation $\mtx{T}$ to $T$ via
$$
\mtx{T} = \mtx{S}\,\left(\tfrac{1}{2}\mtx{I} + \mtx{D}^{*}\right)^{-1}.
$$

For this example, we evaluated an additional error metric by testing the computed NtD operator
to an exact solution $u_{\rm exact}$ to (\ref{eq:BVP}). The function $u_{\rm exact}$ is given
as the potential from a collection of five randomly placed charges inside $\Gamma$, and then
$\vct{f}_{\rm exact}$ and $\vct{g}_{\rm exact}$ are simply the evaluations of $u_{\rm exact}$
and its normal derivative on the quadrature nodes on $\Gamma$. Then the new error measure is
given by
$$
E_{\rm pot} = \frac{\|\vct{f}_{\rm exact} - \mtx{T}_{\rm approx}\,\vct{g}_{\rm exact}\|_{\rm max}}
                   {\|\vct{f}_{\rm exact}\|_{\rm max}},
$$
where $\|\cdot\|_{\rm max}$ is the maximum norm, and where $\mtx{T}_{\rm approx}$ is the
compressed representation of $T$ determined by the randomized sampling scheme proposed.

\begin{figure}
\includegraphics[width=\textwidth]{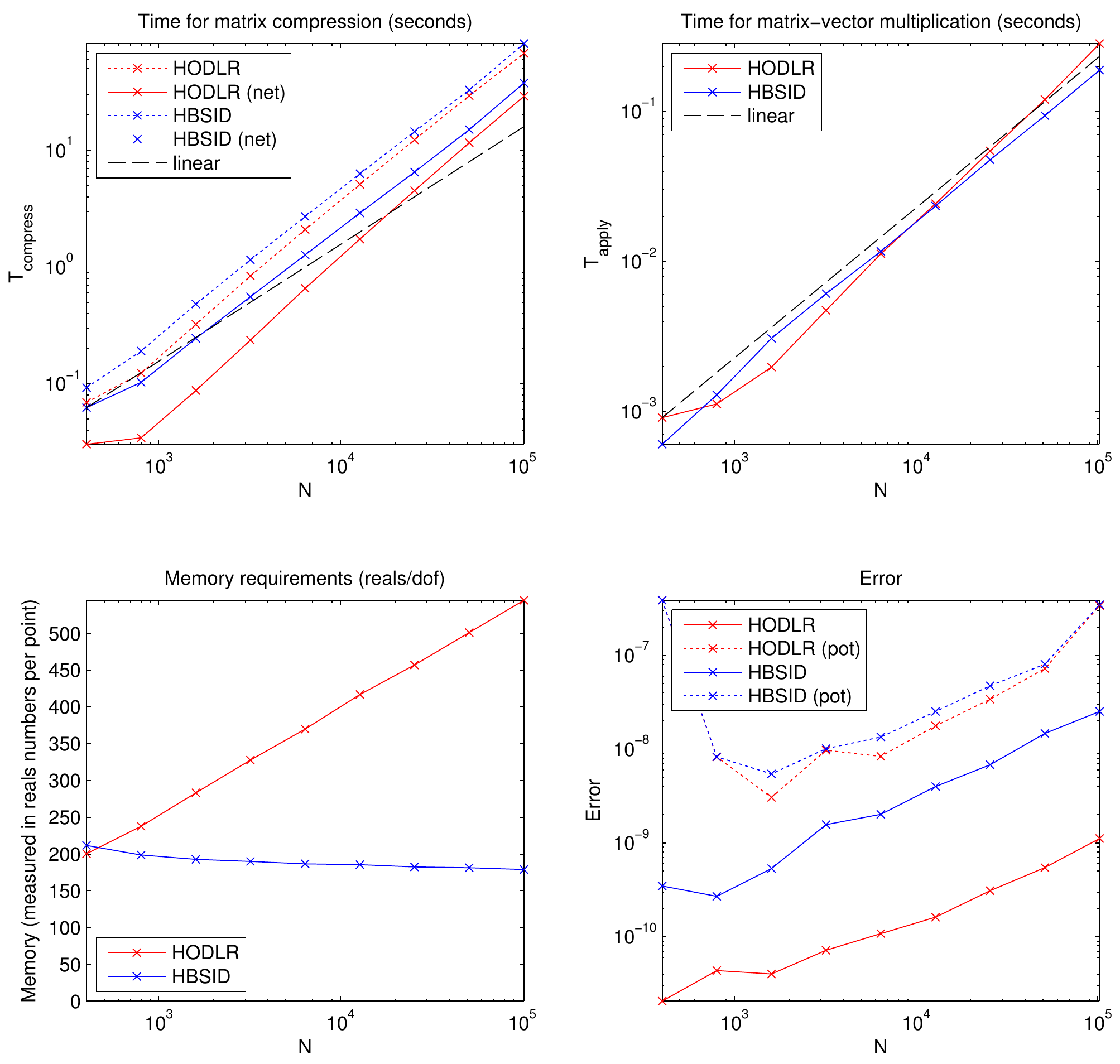}
\caption{Visualization of the results presented in Tables \ref{table:DfN_HODLR}, and \ref{table:DfN_HBSID},
pertaining to the example in Section \ref{sec:DfN}.}
\label{fig:DfN_results}
\end{figure}

The numerical results are
presented in tables \ref{table:DfN_HODLR} and \ref{table:DfN_HBSID},
and summarized in Figure \ref{fig:DfN_results}.

\begin{remark}
\label{remark:DfN}
The formula (\ref{eq:DfN}) for the NtD operator is derived as follows:
We first look for a solution to (\ref{eq:BVP}) of the form $u = Sq$.
Then it can be shown that $q$ must satisfy $(1/2)q + D^{*}q = g$.
Solving for $q$ and using $f = Sq$, we obtain (\ref{eq:DfN}).
\end{remark}

\begin{table}
\begin{center}
\begin{tabular}{|r|r|r|r|r|r|r|r|r|r|r|}\hline
$N$ & $N_{\rm matvec}$ & $T_{\rm compress}$ & $T_{\rm net}$ & $T_{\rm app}$ & $M$ & $M/N$ & $E$ & $E_{\rm pot}$ & $k$\\
&& (sec)& (sec) & (sec) & (MB) & (reals) & & \\ \hline
   400 &   2 x 60 &     0.069 &      0.030 &     0.001 &    0.6 & 200.3 &  2.06e-11  &  3.87e-07 &  32 \\
   800 &   3 x 60 &     0.123 &      0.034 &     0.001 &    1.5 & 237.7 &  4.35e-11  &  8.22e-09 &  35 \\
  1600 &   4 x 60 &     0.323 &      0.088 &     0.002 &    3.5 & 283.1 &  4.00e-11  &  3.06e-09 &  40 \\
  3200 &   5 x 60 &     0.839 &      0.236 &     0.005 &    8.0 & 327.7 &  7.18e-11  &  9.74e-09 &  44 \\
  6400 &   6 x 60 &     2.093 &      0.657 &     0.011 &   18.1 & 370.0 &  1.08e-10  &  8.37e-09 &  45 \\
 12800 &   7 x 60 &     5.108 &      1.734 &     0.024 &   40.7 & 417.0 &  1.61e-10  &  1.76e-08 &  49 \\
 25600 &   8 x 60 &    12.302 &      4.516 &     0.055 &   89.3 & 457.2 &  3.09e-10  &  3.41e-08 &  52 \\
 51200 &   9 x 60 &    29.380 &     11.636 &     0.120 &  195.8 & 501.2 &  5.44e-10  &  7.22e-08 &  54 \\
102400 &  10 x 60 &    67.879 &     29.072 &     0.284 &  426.0 & 545.3 &  1.11e-09  &  3.39e-07 &  55 \\ \hline
\end{tabular}
\end{center}
\caption{Compression to the HODLR format using Algorithm \ref{alg:HODLR} of the NtD operator described in
Section \ref{sec:DfN}. Here $\varepsilon = 10^{-9}$ and $\ell = 60$.}
\label{table:DfN_HODLR}
\end{table}

\begin{table}
\begin{center}
\begin{tabular}{|r|r|r|r|r|r|r|r|r|r|r|}\hline
$N$ & $N_{\rm matvec}$ & $T_{\rm compress}$ & $T_{\rm net}$ & $T_{\rm app}$ & $M$ & $M/N$ & $E$ & $E_{\rm pot}$ & $k$\\
&& (sec)& (sec) & (sec) & (MB) & (reals) & & \\ \hline
   400 &   2 x 60 &     0.093 &      0.062 &     0.001 &    0.6 & 211.7 &  3.48e-10 &  3.86e-07 &  29 \\
   800 &   3 x 60 &     0.190 &      0.103 &     0.001 &    1.2 & 198.7 &  2.70e-10 &  8.25e-09 &  31 \\
  1600 &   4 x 60 &     0.482 &      0.245 &     0.003 &    2.4 & 192.7 &  5.32e-10 &  5.42e-09 &  33 \\
  3200 &   5 x 60 &     1.154 &      0.557 &     0.006 &    4.6 & 189.9 &  1.57e-09 &  1.01e-08 &  36 \\
  6400 &   6 x 60 &     2.708 &      1.269 &     0.012 &    9.1 & 186.5 &  2.02e-09 &  1.35e-08 &  37 \\
 12800 &   7 x 60 &     6.310 &      2.906 &     0.023 &   18.1 & 185.4 &  3.98e-09 &  2.51e-08 &  40 \\
 25600 &   8 x 60 &    14.495 &      6.511 &     0.048 &   35.6 & 182.3 &  6.80e-09 &  4.72e-08 &  41 \\
 51200 &   9 x 60 &    32.906 &     15.033 &     0.094 &   70.8 & 181.4 &  1.47e-08 &  8.04e-08 &  42 \\
102400 &  10 x 60 &    82.238 &     37.772 &     0.189 &  139.8 & 178.9 &  2.51e-08 &  3.46e-07 &  44 \\ \hline
\end{tabular}
\end{center}
\caption{Compression to the HODLR format using Algorithm \ref{alg:HODLR} of the NtD operator described in
Section \ref{sec:DfN}. Here $\varepsilon = 10^{-9}$ and $\ell = 60$.}
\label{table:DfN_HBSID}
\end{table}

\subsection{Dimensional reduction in the Fast Multipole Method}
\label{sec:FMM}
\begin{figure}
\includegraphics[width=\textwidth]{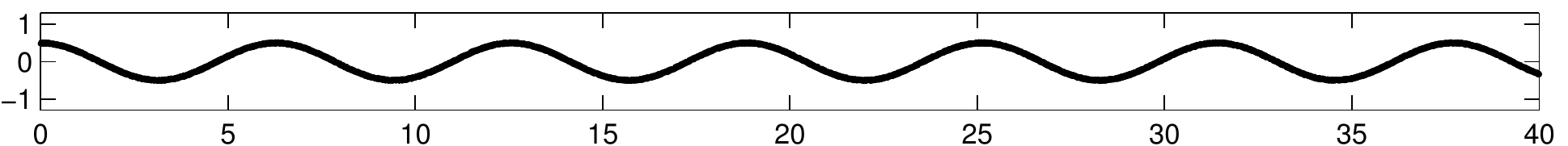}\\
\includegraphics[width=\textwidth]{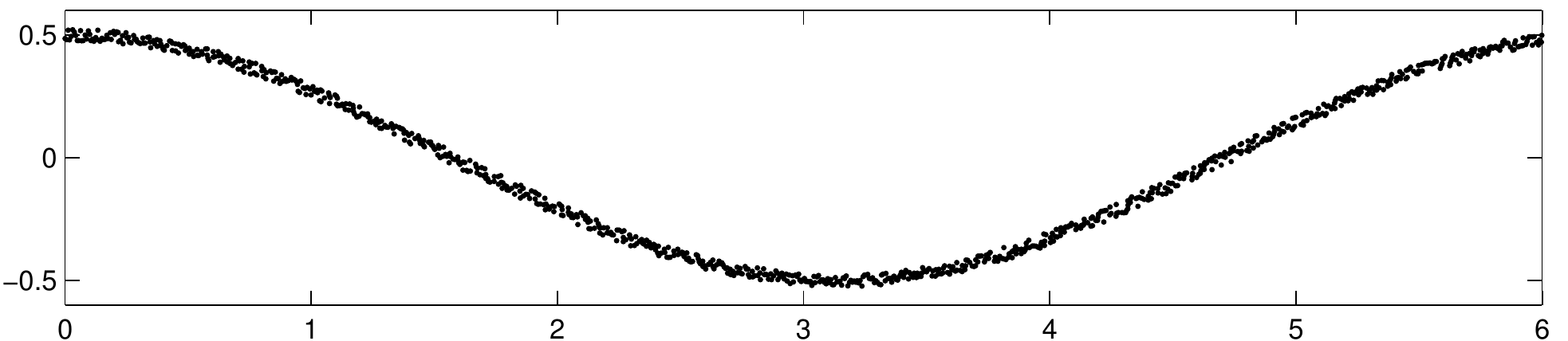}
\caption{Geometry of the potential evaluation map in Section \ref{sec:FMM}.
The top figure shows the entire geometry, and the lower figure shows
a magnification of a small part.}
\label{fig:FMM_geom}
\end{figure}

In this section, we investigate a numerical example where $\mtx{A}$ is a potential evaluation
map for a set of electric charges in the plane. To be precise, we let $\{\pvct{x}_{i}\}_{i=1}^{N}$
denote a set of points located as shown in Figure \ref{fig:FMM_geom}. Then $\mtx{A}$ is the
$N\times N$ matrix with entries
$$
\mtx{A}(i,i) =
\left\{\begin{array}{ll}
\displaystyle -\frac{1}{2\pi}\log|\pvct{x}_{i} - \pvct{x}_{j}|,\qquad&\mbox{when}\ i \neq j,\\
\displaystyle 0&\mbox{when}\ i=j.
\end{array}\right.
$$
We apply $\mtx{A}$ rapidly using the classical Fast Multipole Method \cite{rokhlin1987} with
30th order expansions, which ensures that the error in the black-box code is far smaller than
our requested precision of $\varepsilon = 10^{-9}$. Our FMM is implemented in Matlab as described
in \cite{2011_martinsson_randomhudson}. This implementation is quite inefficient, but has the
ability to apply $\mtx{A}$ to an $N\times \ell$ matrix rather than single vector.

For this example, we use an additional error measure $E_{\rm pot}$ that compares the
compressed matrix $\mtx{A}_{\rm compress}$ to the exact matrix $\mtx{A}$,
evaluated at a subset of the target points. To be precise, we picked at random a subset
of 10 target locations, marked by the index vector $I \subset \{1,\,2,\,\dots,\,N\}$.
We then drew a sequence of ten vectors $\{\vct{\omega}_{i}\}_{i=1}^{10}$ in which each
entry is drawn at random from a uniform distribution on the interval $[0,1]$ (observe
that every charge is \textit{positive}). Then $E_{\rm pot}$ is defined via
$$
E_{\rm pot} = \max_{1 \leq i \leq 10}
\frac{||\mtx{A}(I,:)\vct{\omega}_{i} - \mtx{A}_{\rm compressed}(I,:)\vct{\omega}_{i}||}
     {||\mtx{A}(I,:)\vct{\omega}_{i}||}.
$$

\begin{figure}
\includegraphics[width=\textwidth]{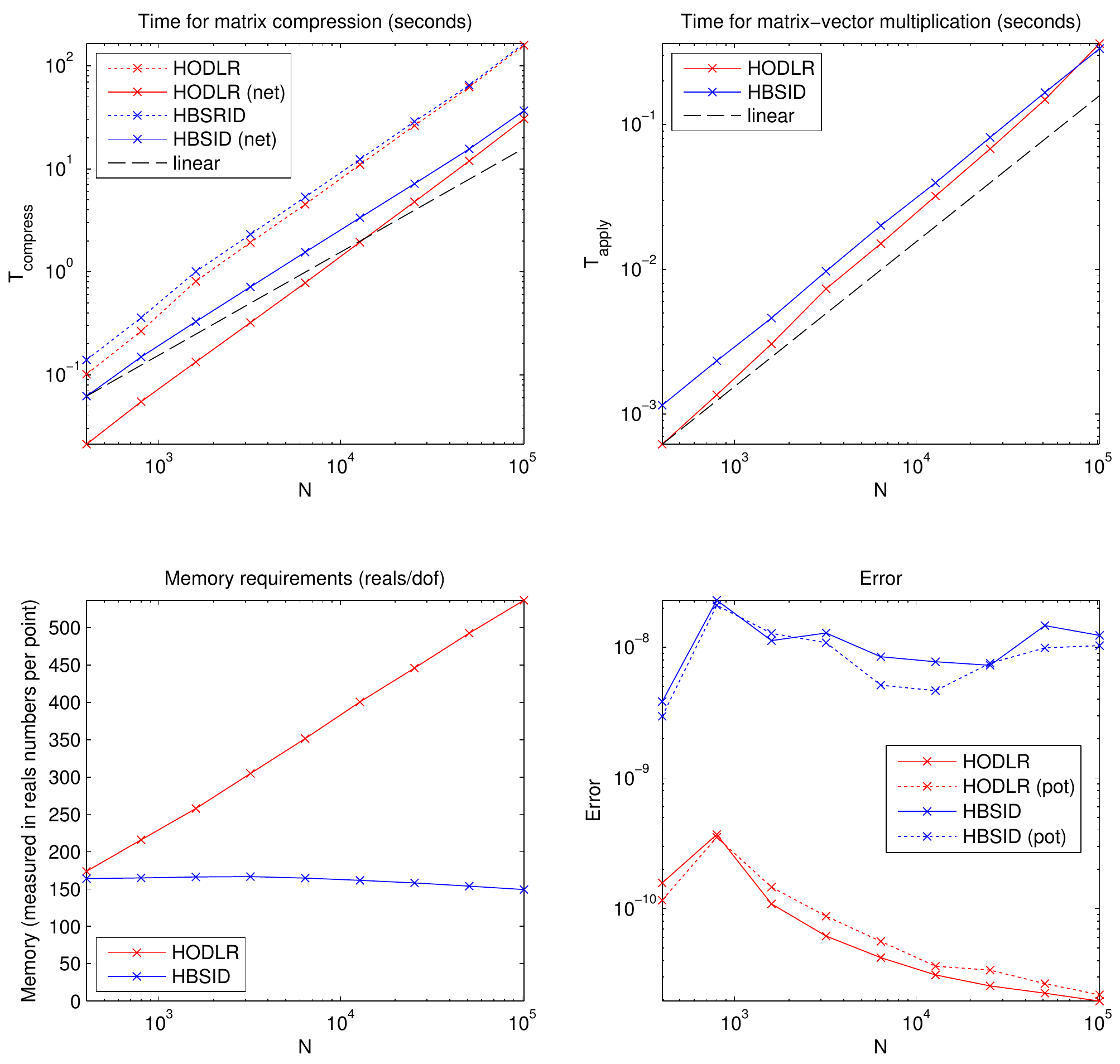}
\caption{Visualization of the results presented in Tables \ref{table:FMM_HODLR}, and \ref{table:FMM_HBSID},
pertaining to the example in Section \ref{sec:FMM}.}
\label{fig:FMM_results}
\end{figure}

The numerical results are
presented in tables \ref{table:FMM_HODLR} and \ref{table:FMM_HBSID},
and summarized in Figure \ref{fig:FMM_results}.

\begin{table}
\begin{center}
\begin{tabular}{|r|r|r|r|r|r|r|r|r|r|r|}\hline
$N$ & $N_{\rm matvec}$ & $T_{\rm compress}$ & $T_{\rm net}$ & $T_{\rm app}$ & $M$ & $M/N$ & $E$ & $E_{\rm pot}$ & $k$\\
&& (sec)& (sec) & (sec) & (MB) & (reals) & & \\ \hline
   400 &   3 x 45 &     0.101 &      0.021 &     0.001 &    0.5 & 173.9 &  1.57e-10  &  1.16e-10 &  21 \\
   800 &   4 x 45 &     0.265 &      0.055 &     0.001 &    1.3 & 215.9 &  3.69e-10  &  3.52e-10 &  23 \\
  1600 &   5 x 45 &     0.812 &      0.133 &     0.003 &    3.1 & 257.8 &  1.09e-10  &  1.46e-10 &  24 \\
  3200 &   6 x 45 &     1.925 &      0.321 &     0.007 &    7.4 & 304.9 &  6.17e-11  &  8.77e-11 &  26 \\
  6400 &   7 x 45 &     4.509 &      0.783 &     0.015 &   17.2 & 351.5 &  4.22e-11  &  5.61e-11 &  28 \\
 12800 &   8 x 45 &    10.980 &      1.941 &     0.032 &   39.1 & 400.6 &  3.11e-11  &  3.63e-11 &  30 \\
 25600 &   9 x 45 &    26.271 &      4.788 &     0.068 &   87.1 & 446.0 &  2.56e-11  &  3.38e-11 &  30 \\
 51200 &  10 x 45 &    62.119 &     11.938 &     0.149 &  192.5 & 492.7 &  2.26e-11  &  2.67e-11 &  32 \\
102400 &  11 x 45 &   158.617 &     30.621 &     0.362 &  419.5 & 537.0 &  1.96e-11  &  2.19e-11 &  32 \\ \hline
\end{tabular}
\end{center}
\caption{Compression to the HODLR format using Algorithm \ref{alg:HODLR} of the potential evaluation matrix
described in Section \ref{sec:FMM}. Here $\varepsilon = 10^{-9}$ and $\ell = 45$.}
\label{table:FMM_HODLR}
\end{table}

\begin{table}
\begin{center}
\begin{tabular}{|r|r|r|r|r|r|r|r|r|r|r|}\hline
$N$ & $N_{\rm matvec}$ & $T_{\rm compress}$ & $T_{\rm net}$ & $T_{\rm app}$ & $M$ & $M/N$ & $E$ & $E_{\rm pot}$ & $k$\\
&& (sec)& (sec) & (sec) & (MB) & (reals) & & \\ \hline
   400 &   3 x 45 &     0.140 &      0.062 &     0.001 &    0.5 & 164.1 &  3.84e-09 &  2.96e-09 &  31 \\
   800 &   4 x 45 &     0.358 &      0.149 &     0.002 &    1.0 & 164.8 &  2.29e-08 &  2.10e-08 &  33 \\
  1600 &   5 x 45 &     1.007 &      0.329 &     0.005 &    2.0 & 166.1 &  1.12e-08 &  1.28e-08 &  37 \\
  3200 &   6 x 45 &     2.306 &      0.713 &     0.010 &    4.1 & 166.5 &  1.29e-08 &  1.08e-08 &  39 \\
  6400 &   7 x 45 &     5.334 &      1.554 &     0.020 &    8.0 & 164.7 &  8.47e-09 &  5.13e-09 &  40 \\
 12800 &   8 x 45 &    12.366 &      3.350 &     0.040 &   15.8 & 161.7 &  7.75e-09 &  4.66e-09 &  42 \\
 25600 &   9 x 45 &    28.813 &      7.170 &     0.081 &   30.9 & 158.1 &  7.28e-09 &  7.56e-09 &  43 \\
 51200 &  10 x 45 &    64.838 &     15.572 &     0.166 &   60.1 & 154.0 &  1.47e-08 &  9.91e-09 &  44 \\
102400 &  11 x 45 &   164.750 &     36.686 &     0.335 &  116.7 & 149.4 &  1.23e-08 &  1.03e-08 &  43 \\ \hline
\end{tabular}
\end{center}
\caption{Compression to the HBSID format using Algorithm \ref{alg:HBSID} of the potential evaluation matrix
described in Section \ref{sec:FMM}. Here $\varepsilon = 10^{-9}$ and $\ell = 45$.}
\label{table:FMM_HBSID}
\end{table}

\subsection{Compression of frontal matrices in nested dissection}
\label{sec:nesteddiss}
\begin{figure}
\setlength{\unitlength}{1mm}
\begin{picture}(125,35)
\put(003,00){\includegraphics[width=120mm]{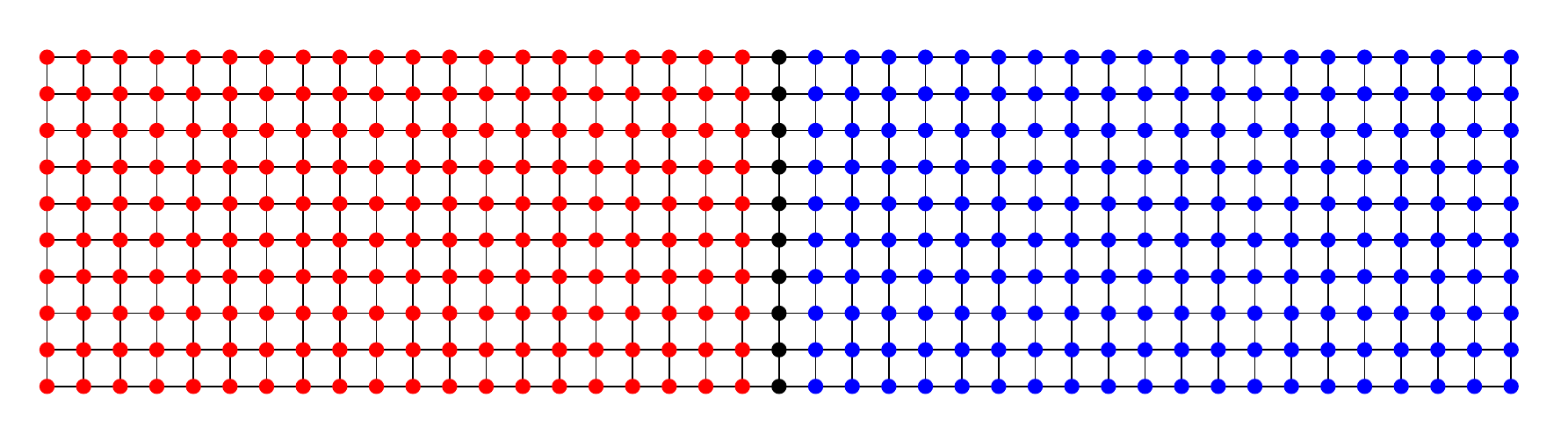}}
\put(000,18){$\color{red}I_{1}$}
\put(122,18){$\color{blue}I_{2}$}
\put(062,31){$I_{3}$}
\end{picture}
\caption{Geometry of problem described in Section \ref{sec:nesteddiss}. We consider a grid
conduction problem on the grid shown. As $N$ is increased, the \textit{width} or the grid
is fixed at 41 nodes, while the \textit{height} of the grid equals $N$.}
\label{fig:nesteddiss_geom}
\end{figure}

Our final example applies the proposed compression schemes to the problem of
constructing $O(N)$ direct solvers for the sparse linear systems arising upon
the discretization of elliptic PDEs via finite difference or finite element
methods. The idea is to build a solver on the classical ``nested dissection''
scheme of George \cite{george_1973,1989_directbook_duff,2006_davis_directsolverbook}.
In standard implementations, the problem of this direct solver is that it requires
the inversion or LU-factorization of a set of successively larger dense matrices.
However, it has recently been demonstrated that while these matrices are dense,
they have internal structure that allows for linear or close to linear time matrix
algebra to be executed, \cite{2007_leborne_HLU,2009_xia_superfast,2009_martinsson_FEM,2013_martinsson_DtN_linearcomplexity}.
In this manuscript, we test the proposed compression scheme on a set of matrices
whose behavior is directly analogous to the matrices encountered in the algorithms
of \cite{2007_leborne_HLU,2009_xia_superfast,2009_martinsson_FEM,2013_martinsson_DtN_linearcomplexity}.
To be precise, let $\mtx{B}$ denote the sparse coefficient matrix associated with a grid conduction
problem on the grid shown in  Figure \ref{fig:nesteddiss_geom}. Each bar has a conductivity that
is drawn at random from a uniform distribution on the interval $[1,2]$.
Let $I_{1},I_{2},I_{3}$ denote three index vectors that mark the three regions shown in
Figure \ref{fig:nesteddiss_geom}, set
$$
\mtx{B}_{ij} = \mtx{B}(I_{i},I_{j}),
\qquad i,j = 1,2,3,
$$
and then define the $N\times N$ matrix $\mtx{A}$ via
\begin{equation}
\label{eq:defA_nesteddiss}
\mtx{A} = \mtx{B}_{33} - \mtx{B}_{31}\mtx{B}_{11}^{-1}\mtx{B}_{13} - \mtx{B}_{32}\mtx{B}_{22}^{-1}\mtx{B}_{23}.
\end{equation}
The relevance of the matrix $\mtx{A}$ is discussed in some detail in Remark \ref{remark:nesteddiss}.

In our numerical experiments, the black-box application of $\mtx{A}$, as defined by (\ref{eq:defA_nesteddiss})
was executed using the sparse matrix built-in routines in Matlab, which relies on UMFPACK \cite{2004_davis_UMFPACK}
for the sparse solves implicit in the application of $\mtx{B}_{11}^{-1}$ and $\mtx{B}_{22}^{-1}$.

\begin{figure}
\includegraphics[width=\textwidth]{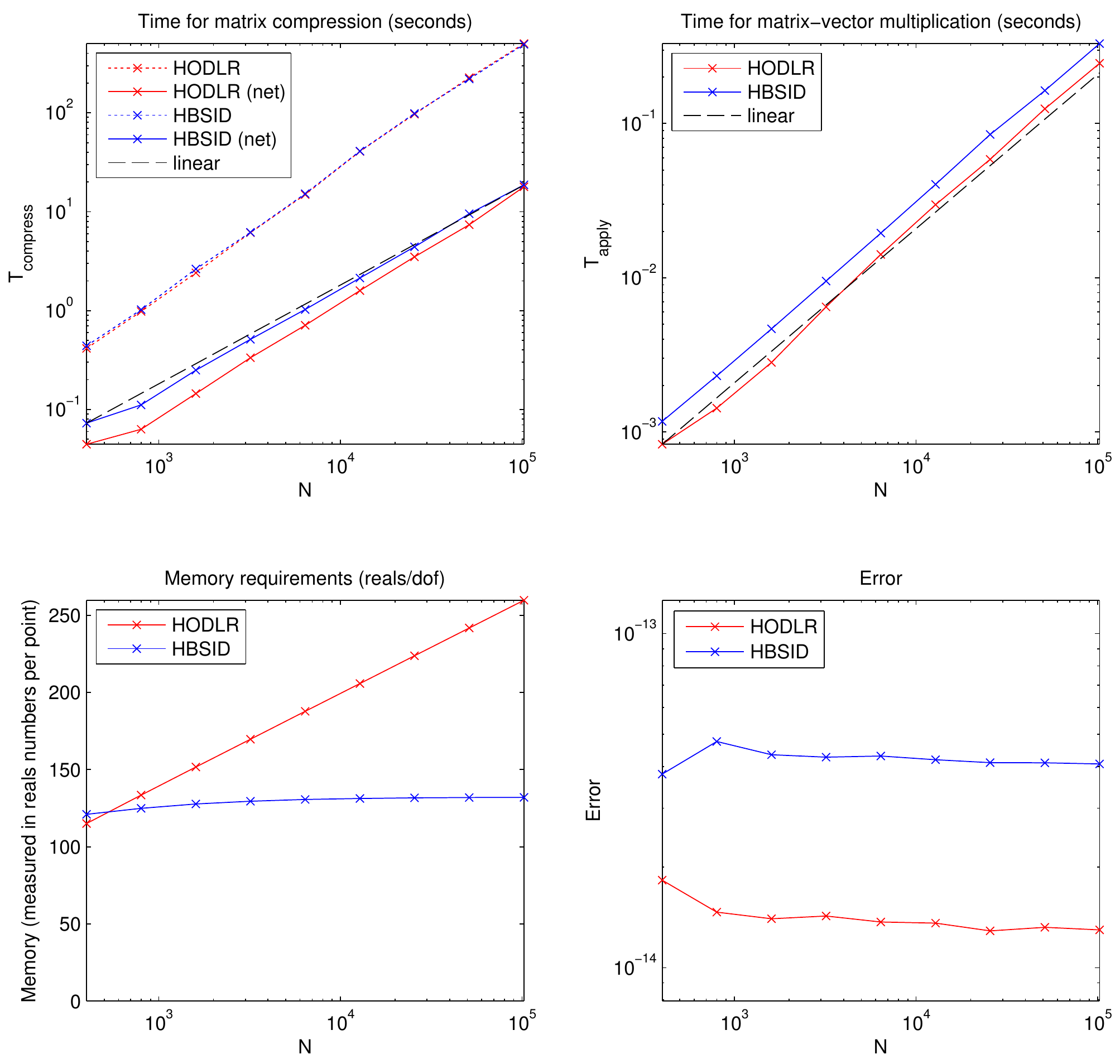}
\caption{Visualization of the results presented in Tables \ref{table:nesteddiss_HODLR}, and \ref{table:nesteddiss_HBSID},
pertaining to the example in Section \ref{sec:nesteddiss}.}
\label{fig:nesteddiss_results}
\end{figure}

The numerical results are
presented in tables \ref{table:nesteddiss_HODLR} and \ref{table:nesteddiss_HBSID},
and summarized in Figure \ref{fig:nesteddiss_results}.

\begin{remark}
\label{remark:nesteddiss}
To illustrate the connection between the matrix $\mtx{A}$, as defined by (\ref{eq:defA_nesteddiss}),
and the LU-factorization of a matrix such as $\mtx{B}$, observe first that the blocks
$\mtx{B}_{12}$ and $\mtx{B}_{21}$ are zero, so that (up to a permutation of the rows and columns),
$$
\mtx{B} =
\left[\begin{array}{ccc}
\mtx{B}_{11} & \mtx{0}      & \mtx{B}_{13} \\
\mtx{0}      & \mtx{B}_{22} & \mtx{B}_{23} \\
\mtx{B}_{31} & \mtx{B}_{32} & \mtx{B}_{33}
\end{array}\right].
$$
Next suppose that we can somehow determine the LU-factorizations of $\mtx{B}_{11}$ and $\mtx{B}_{22}$,
$$
\mtx{B}_{11} = \mtx{L}_{11}\mtx{U}_{11},
\qquad\mbox{and}\qquad
\mtx{B}_{22} = \mtx{L}_{22}\mtx{U}_{22}.
$$
Then the LU-factorization of $\mtx{B}$ is given by
$$
\mtx{B} =
\left[\begin{array}{ccc}
\mtx{L}_{11} & \mtx{0}      & \mtx{0} \\
\mtx{0}      & \mtx{L}_{22} & \mtx{0} \\
\mtx{B}_{31}\mtx{U}_{11}^{-1} & \mtx{B}_{32}\mtx{U}_{22}^{-1} & \mtx{L}_{33}
\end{array}\right]\,
\left[\begin{array}{ccc}
\mtx{U}_{11} & \mtx{0}      & \mtx{L}_{11}^{-1}\mtx{B}_{13} \\
\mtx{0}      & \mtx{U}_{22} & \mtx{L}_{22}^{-1}\mtx{B}_{23} \\
\mtx{0}      & \mtx{0}      & \mtx{U}_{33}
\end{array}\right],
$$
where $\mtx{L}_{33}$ and $\mtx{U}_{33}$ are defined as the LU factors of
\begin{equation}
\label{eq:library}
\mtx{L}_{33}\mtx{U}_{33} = \underbrace{\mtx{B}_{33} - \mtx{B}_{31}\mtx{U}_{11}^{-1}\mtx{L}_{11}^{-1}\mtx{B}_{13} - \mtx{B}_{32}\mtx{U}_{22}^{-1}\mtx{L}_{22}^{-1}\mtx{B}_{23}}_{=:\mtx{A}}.
\end{equation}
Observe that since the matrices $\mtx{L}_{11},\,\mtx{U}_{11},\,\mtx{L}_{22},\,\mtx{U}_{22}$
are all triangular, their inverses are inexpensive to apply. To summarize, \text{if} one
can cheaply evaluate the LU factorization (\ref{eq:library}), then the task of LU-factoring
$\mtx{B}$ directly reduces to the task of LU factoring the two matrices $\mtx{B}_{11}$
and $\mtx{B}_{22}$, which both involve about half as many variables as $\mtx{B}$. The
classical nested dissection idea is now to apply this observation recursively to factor
$\mtx{B}_{11}$ and $\mtx{B}_{22}$. The problem of this scheme has traditionally been that
in order to evaluate $\mtx{L}_{33}$ and $\mtx{U}_{33}$ in (\ref{eq:library}), one must
factorize the \textit{dense} matrix $\mtx{A}$.
\end{remark}

\begin{table}
\begin{center}
\begin{tabular}{|r|r|r|r|r|r|r|r|r|r|r|}\hline
$N$ & $N_{\rm matvec}$ & $T_{\rm compress}$ & $T_{\rm net}$ & $T_{\rm app}$ & $M$ & $M/N$ & $E$ & $k$\\
&& (sec)& (sec) & (sec) & (MB) & (reals) & \\ \hline
   400 &   3 x 25 &     0.348 &      0.012 &     0.001 &    0.4 & 115.0 &  1.83e-14  &   9 \\
   800 &   4 x 25 &     0.849 &      0.025 &     0.001 &    0.8 & 133.4 &  1.47e-14  &   9 \\
  1600 &   5 x 25 &     2.143 &      0.064 &     0.002 &    1.9 & 151.6 &  1.40e-14  &   9 \\
  3200 &   6 x 25 &     5.588 &      0.141 &     0.005 &    4.1 & 169.7 &  1.43e-14  &   9 \\
  6400 &   7 x 25 &    13.926 &      0.336 &     0.011 &    9.2 & 187.8 &  1.37e-14  &   9 \\
 12800 &   8 x 25 &    38.558 &      0.903 &     0.023 &   20.1 & 205.8 &  1.36e-14  &   9 \\
 25600 &   9 x 25 &    91.382 &      2.322 &     0.048 &   43.7 & 223.8 &  1.29e-14  &   9 \\
 51200 &  10 x 25 &   208.704 &      5.874 &     0.098 &   94.5 & 241.8 &  1.32e-14  &   9 \\
102400 &  11 x 25 &   468.981 &     14.203 &     0.204 &  203.0 & 259.8 &  1.30e-14  &   9 \\ \hline
\end{tabular}
\end{center}
\caption{Compression to the HODLR format using Algorithm \ref{alg:HODLR} of a simulated ``frontal
matrix'' in the nested dissection technique, as described in Section \ref{sec:nesteddiss}. Here $\varepsilon = 10^{-9}$ and $\ell = 25$.}
\label{table:nesteddiss_HODLR}
\end{table}

\begin{table}
\begin{center}
\begin{tabular}{|r|r|r|r|r|r|r|r|r|r|r|}\hline
$N$ & $N_{\rm matvec}$ & $T_{\rm compress}$ & $T_{\rm net}$ & $T_{\rm app}$ & $M$ & $M/N$ & $E$ & $k$\\
&& (sec)& (sec) & (sec) & (MB) & (reals) & \\ \hline
   400 &   3 x 25 &     0.361 &      0.032 &     0.001 &    0.4 & 121.0 &  3.80e-14 &  18 \\
   800 &   4 x 25 &     0.867 &      0.066 &     0.002 &    0.8 & 124.9 &  4.75e-14 &  18 \\
  1600 &   5 x 25 &     2.180 &      0.150 &     0.005 &    1.6 & 127.7 &  4.34e-14 &  18 \\
  3200 &   6 x 25 &     5.575 &      0.322 &     0.010 &    3.2 & 129.5 &  4.27e-14 &  18 \\
  6400 &   7 x 25 &    14.470 &      0.708 &     0.019 &    6.4 & 130.6 &  4.30e-14 &  18 \\
 12800 &   8 x 25 &    38.945 &      1.525 &     0.040 &   12.8 & 131.3 &  4.19e-14 &  18 \\
 25600 &   9 x 25 &    91.563 &      3.253 &     0.081 &   25.7 & 131.7 &  4.11e-14 &  18 \\
 51200 &  10 x 25 &   208.145 &      7.440 &     0.157 &   51.5 & 131.9 &  4.10e-14 &  18 \\
102400 &  11 x 25 &   465.309 &     15.759 &     0.314 &  103.1 & 132.0 &  4.07e-14 &  18 \\ \hline
\end{tabular}
\end{center}
\caption{Compression to the HBSID format using Algorithm \ref{alg:HBSID} of a simulated ``frontal
matrix'' in the nested dissection technique, as described in Section \ref{sec:nesteddiss}.
Here $\varepsilon = 10^{-9}$ and $\ell = 25$.}
\label{table:nesteddiss_HBSID}
\end{table}

\subsection{Summary of observations from numerical experiments}
To close this section, we make some observations and conjectures:
\begin{itemize}
\item In all examples examined, numerical evidence supports the claims on asymptotic scaling
made in Sections \ref{sec:HODLRasympt} and \ref{sec:HBSasympt}.
\item Excellent approximation errors are obtained in every case.
Aggregation of errors over levels is a very minor problem.
\item The computational time is in all cases dominated by the time
spent in the external black-box multiplier. As a consequence, the
primary route by which the proposed algorithm could be improved
would be to lessen the number of matrix-vector multiplications
required.
\item For modest problem sizes, the HODLR algorithm is very fast
and easy to use. However, as problems grow large, the memory
requirements of the HODLR format become slightly problematic,
and the HBSID algorithm gains a more pronounced advantage.
\end{itemize}

\begin{figure}
\fbox{\begin{minipage}{190mm}
      \begin{tabbing}
      \hspace{5mm} \= \hspace{5mm} \= \hspace{5mm} \= \hspace{5mm}\kill
      \textit{\color{blue}Build outgoing expansions on level $m$.}\\
      \textbf{loop} over all nodes $\tau$ on level $m$\\
      \> $\hat{\vct{q}}_{\tau} = \mtb{V}_{\tau}^{*}\,\vct{q}(I_{\tau})$\\
      \textbf{end loop} \\[2mm]
      \textit{\color{blue}Build outgoing expansions on levels coarser than $m$ (upwards pass).}\\
      \textbf{loop} over levels $\ell = (m-1):(-1):1$\\
      \> \textbf{loop} over boxes $\tau$ on level $\ell$\\
      \> \> Let $\{\alpha,\beta\}$ denote the children of $\tau$.\\
      \> \> $\hat{\vct{q}}_{\tau} = \mtx{V}_{\tau}^{*}
                                    \left[\begin{array}{r}
                                    \hat{\vct{q}}_{\alpha} \\ \hat{\vct{q}}_{\beta}
                                    \end{array}\right]$.\\
      \> \textbf{end loop} \\
      \textbf{end loop} \\[2mm]
      \textit{\color{blue}Build incoming expansions for the children of the root.}\\
      Let $\{\alpha,\beta\}$ denote the children of the root node.\\
      $\hat{\vct{u}}_{\alpha} = \mtx{B}_{\alpha,\beta}\,\hat{\vct{q}}_{\beta }$.\\
      $\hat{\vct{u}}_{\beta } = \mtx{B}_{\beta,\alpha}\,\hat{\vct{q}}_{\alpha}$.\\[2mm]
      \textit{\color{blue}Build incoming expansions on levels coarser than $m$ (downwards pass).}\\
      \textbf{loop} over levels $\ell = (m-1):(-1):1$\\
      \> \textbf{loop} over boxes $\tau$ on level $\ell$\\
      \> \> Let $\{\alpha,\beta\}$ denote the children of $\tau$.\\
      \> \> $\hat{\vct{u}}_{\alpha} = \mtx{B}_{\alpha,\beta}\,\hat{\vct{q}}_{\beta } +
                                      \mtx{U}_{\tau}(J_{\alpha},:)\,\hat{\vct{u}}_{\tau}$.\\
      \> \> $\hat{\vct{u}}_{\beta } = \mtx{B}_{\beta,\alpha}\,\hat{\vct{q}}_{\alpha} +
                                      \mtx{U}_{\tau}(J_{\beta },:)\,\hat{\vct{u}}_{\tau}$.\\
      \> \textbf{end loop} \\
      \textbf{end loop} \\[2mm]
      \textit{\color{blue}Build incoming expansions on level $m$.}\\
      \textbf{loop} over boxes $\tau$ on level $m$\\
      \> $\vct{u}(I_{\tau}) = \mtb{U}_{\tau}\,\hat{\vct{u}}_{\tau}$\\
      \textbf{end loop}
      \end{tabbing}
      \end{minipage}}
\caption{Application of $\mtx{A}^{(m)}$ in the HBS framework. Given a vector $\vct{q}$
of charges, build the vector $\vct{u} = \mtx{A}^{(m)}\,\vct{q}$ of potentials.}
\label{alg:apply_A_partial}
\end{figure}

\bibliographystyle{siam}
\bibliography{main_bib}

\end{document}